\documentclass[12pt]{amsart}
\usepackage{subfig}
\usepackage{hyperref}
\usepackage{enumitem}
\usepackage{graphicx}
\usepackage{mathrsfs}
\usepackage{xfrac}
\usepackage[T1]{fontenc}
\usepackage{amssymb}
\usepackage[cmtip,all]{xy}
\usepackage{tikz-cd}
\usepackage{tikzsymbols}
\usepackage{amscd}
\usepackage{mathtools}
\usepackage[utf8]{inputenc}
\usepackage{color} 
\usepackage[centertags]{amsmath}
\usepackage{indentfirst}
\usepackage[all]{xy}
\usepackage{multicol}
\usepackage{upgreek}
\usepackage{tcolorbox}
\usepackage{relsize,fancyvrb}
\usetikzlibrary{shapes,snakes}
\usepackage{relsize}
\usepackage{graphicx}
\usepackage{tikz}
\usepackage{stmaryrd}
\usepackage{xcolor}
\usepackage{pdflscape}
\usepackage{tikz}
\usepackage{scalerel}[2016/12/29]
\usetikzlibrary{arrows,decorations.pathmorphing}

\makeatletter
\@namedef{subjclassname@2020}{%
	\textup{2020} Mathematics Subject Classification}
\makeatother

\newcommand\ka{{$\boldsymbol{\left( \Omega\right) }$}}

\newtcolorbox{mybox}{colback=red!5!white,colframe=red!50!black} 
\newtcbox{\myboxx}{on line,
	colframe=red!75!black,colback=red!5!white,
	boxrule=1pt,arc=4pt,boxsep=0pt,left=6pt,right=6pt,top=6pt,bottom=6pt}
\newtcolorbox{myboxxx}{colback=red!4!white,colframe=red!75!black}

\newtheorem{thm}{Theorem}[section]
\newtheorem{defnt}[thm]{Definition}
\newtheorem{prop}[thm]{Proposition}

\newtheorem{cor}[thm]{Corollary}
\newtheorem{lemma}[thm]{Lemma}
\newtheorem{question}[thm]{Question}

\newtheorem{remark}[thm]{Remark}
\newenvironment{proof*}{\begin{trivlist}\item[\hskip\labelsep{\em Proof of Theorem \ref{t2}}]}{%
		\hfill$\square$\rm\end{trivlist}}
\newenvironment{proof**}{\begin{trivlist}\item[\hskip\labelsep{\em Proof of Proposition \ref{p4}}]}{%
		\hfill$\square$\rm\end{trivlist}}
\newenvironment{proof***}{\begin{trivlist}\item[\hskip\labelsep{\em Proof of Proposition \ref{p41}}]}{%
		\hfill$\square$\rm\end{trivlist}}
\newenvironment{proof****}{\begin{trivlist}\item[\hskip\labelsep{\em Proof of Propostion \ref{p42}}]}{%
		\hfill$\square$\rm\end{trivlist}}

\definecolor{forestgreen}{rgb}{0.13, 0.55, 0.13} 
\definecolor{frenchblue}{rgb}{0.0, 0.45, 0.73}
\definecolor{amber}{rgb}{1.0, 0.75, 0.0}
\colorlet{Mycolor1}{green!20!orange!50!purple!50}
\colorlet{Mycolor2}{red!75!black}
\colorlet{Mycolor4}{blue!75!black}

\tikzset{My Style1/.style={black, draw=blue, fill=blue!20, minimum size=1.75cm}}
\tikzset{My Style2/.style={black, draw=red, fill=red!20, minimum size=1.75cm}}
\tikzset{My Style3/.style={black, draw=amber, fill=amber, minimum size=0.55cm}}
\tikzset{My Style5/.style={black, draw=blue, fill=blue!20, minimum size=2cm}}
\tikzset{My Style4/.style={black, draw=red, fill=red!20, minimum size=2cm}}
\tikzset{My Style6/.style={black, draw=black, fill=white, minimum size=1.85cm}}

\tikzset{Mynew Style1/.style={black, draw=blue, fill=blue!20, minimum size=0.55cm}}
\tikzset{Mynew Style2/.style={black, draw=red, fill=red!20, minimum size=0.55cm}}
\tikzset{Mynew Style3/.style={black, draw=black, fill=white, minimum size=0.55cm}}
\tikzset{Mynew Style4/.style={black, draw=amber!75!black, fill=amber!50, minimum size=0.55cm}}

\tikzset{Myneww Style1/.style={black, draw=blue, fill=blue!20, minimum size=1.25cm}}
\tikzset{Myneww Style2/.style={black, draw=red, fill=red!20, minimum size=1.25cm}}
\tikzset{Myneww Style3/.style={black, draw=black, fill=white, minimum size=1.25cm}}

\numberwithin{equation}{section}

\begin{document}
	
	
	\baselineskip=17pt
	
	
	\title[On Power Series Subspaces of Certain Nuclear Fr\'echet Spaces]{On Power Series Subspaces of Certain Nuclear Fr\'echet Spaces}
	
	\author[N. Doğan]{Nazlı Doğan}
	\address{Fatih Sultan Mehmet Vakıf University\\Beyo\u{g}lu, Istanbul, Turkey}
	\email{ndogan@fsm.edu.tr}
	
	\date{}
	
	\begin{abstract}
		The diametral dimension, $\Delta(E)$, and the approximate diametral dimension, $\delta(E)$ of an element $E$ of a large class of nuclear Fr\'echet spaces are set theoretically between the corresponding invariant of power series spaces $\Lambda_{1}(\varepsilon)$ and $\Lambda_{\infty}(\varepsilon)$ for some exponent sequence $\varepsilon$. Aytuna et al., \cite{AKT2}, proved that $E$ contains a complemented subspace which is isomorphic to $\Lambda_{\infty}(\varepsilon)$ provided $\Delta(E)=\Delta( \Lambda_{\infty}(\varepsilon))$ and $\varepsilon$ is stable. In this article, we will consider the other extreme case and we proved that in this large family, there exist nuclear Fr\'echet spaces, even regular nuclear Köthe spaces, satisfying $\Delta(E)=\Delta(\Lambda_{1}(\varepsilon))$ such that there is no subspace of $E$ which is isomorphic to $\Lambda_{1}(\varepsilon)$. 
	\end{abstract}
	\maketitle
	\section{Introduction}
	\par Fr\'echet spaces are one of the leading class of locally convex spaces and include most of the important examples of non-normable locally convex spaces. Power series spaces constitute a well studied class in the theory of Fr\'echet spaces. Subspaces and quotient spaces of a nuclear
	stable power series space are characterized by Vogt and Wagner (\cite{VDN}, \cite{VW}, \cite{VW2}) in terms of diametral dimension and DN-$\Omega$ type linear-topological invariants. The topological invariants \underline{DN} and $\Omega$ are enjoyed by many natural nuclear Fr\'echet
	spaces appearing in analysis and these invariants play an important role in this study.
	\par  Let $E$ be a nuclear Fr\'echet space which satisfies $\underline{DN}$ and $\Omega$. Then it is a known fact that the diametral dimension $\Delta (E)$ and the approximate diametral dimension $\delta (E)$ of $E$ are set theoretically between corresponding invariant of power series spaces $\Lambda_{1}(\varepsilon)$ and $\Lambda_{\infty}\left(\varepsilon\right)$ for some specific exponent sequence $\varepsilon$. The sequence $\varepsilon$ is called  associated exponent sequence of E. In \cite{AKT2}, Aytuna et al. proved that a nuclear Fr\'echet space $E$ with the properties $\underline{DN}$ and $\Omega$ contains a complemented copy of $\Lambda_{\infty}\left(\varepsilon\right)$ provided the diametral dimensions of $E$ and $\Lambda_{\infty}\left(\varepsilon\right)$ are equal and $\varepsilon$ is stable. In this article, we deal with the other extreme, namely, the main question in this article is:
	\begin{question}\label{qi1} Let $E$ be a nuclear Fr\'echet space with the properties $\underline{DN}$ and $\Omega$ and $\varepsilon$ be the associated exponent sequence of $E$. Is there a (complemented) subspace of $E$ which is isomorphic to $\Lambda_{1}(\varepsilon)$ if $\Delta\left(E\right)=\Delta\left(\Lambda_{1}\left( \varepsilon\right)\right)$?
	\end{question}
	\par This problem led us to examine the relationship between the diametral dimension and the other invariants. The most appropriate topological invariants for comparison with the diametral dimension is the approximate diametral dimension. Then, we ask the following question:
	\begin{question}\label{qi2} Let $E$ be a nuclear Fr\'echet space with the properties $\underline{DN}$ and $\Omega$. If diametral dimension of $E$ coincides with that of a power series space, then does this imply that the approximate diametral dimension also do the same and vice versa?
	\end{question}
	\par
	In \cite{ND}, we showed that Question \ref{qi2} has an affirmative answer when power series space is of infinite type. Then we searched an answer for the Question \ref{qi2} in the finite type case and, in this regard, we first proved that the condition $\delta\left(E\right)= \delta\left(\Lambda_{1}\left(\varepsilon\right)\right)$ always implies $\Delta\left(E\right)= \Delta\left(\Lambda_{1}\left(\varepsilon\right)\right)$. We also constructed some sufficient conditions to prove the other direction. It turns out that the existence of a prominent bounded subset in the nuclear Fr\'echet space $E$ plays a decisive role for the answer of Question \ref{qi2}. In \cite[Theorem 4.8]{ND}, we proved that $\delta\left(E\right)=\delta\left(\Lambda_{1}\left(\varepsilon\right)\right)$ if and only if $E$ has a prominent bounded set and $\Delta\left(E\right)=\Delta(\Lambda_{1}\left(\varepsilon\right))$.
	\par In this article, after giving some preliminary materials in Section 2, we construct a family $\bf \mathcal{K}$ of nuclear K\"{o}the spaces $K(a_{k,n})$ parametrized by a sequence $\alpha$ satisfying the properties $\underline{DN}$ and $\Omega$. First we show that for an element of the family of $\bf \mathcal{K}$ which is parameterized by a stable sequence $\alpha$, $\Delta(K(a_{k,n}))=\Delta(\Lambda_{1}(\alpha))$ and $\delta(K(a_{k,n}))= \delta(\Lambda_{1}(\alpha))$. Second, we prove that for any element of the family of $\bf \mathcal{K}$  which is parameterized by an unstable sequence $\alpha$, $\Delta(K(a_{k,n}))=\Delta(\Lambda_{1}(\varepsilon))$ and $\delta(K(a_{k,n}))\neq \delta(\Lambda_{1}(\varepsilon))$ for its associated exponent sequence $\varepsilon$. This show that the Question \ref{qi2} has a negative answer for power series space of finite type. Furthermore, we prove in Theorem \ref{t3} that the first question has a negative answer, that is, $\Lambda_{1}(\varepsilon)$ is not isomorphic to any subspace of these K\"{o}the spaces $K(a_{k,n})$, let alone is isomorphic to a complemented subspace, though the condition $\Delta(K(a_{k,n}))=\Delta(\Lambda_{1}(\varepsilon))$ is satisfied. Motivated by our finding in \cite{ND}, we compile some additional information, for instance, for any element $E$ of the family $\bf \mathcal{K}$ parameterized by an unstable sequence,  
	\begin{itemize}
		\item[1.] $E$ does not have a prominent bounded set.
		\item[2.] Although the equality $\Delta(E)=\Lambda_{1}(\varepsilon)$ is satisfied and the canonical imbedding from $\Delta(E)$ into $\Lambda_{1}(\varepsilon)$ has a closed graph, the canonical imbedding from $\Delta(E)$ into $\Lambda_{1}(\varepsilon)$ is not continuous.
	\end{itemize}
	\section{PRELIMINARIES}
	\par In this section, after establishing terminology and notation, we collect some basic facts and definitions that are needed them in the sequel.
	\par Throughout the article, $E$ will denote a nuclear Fr\'echet space with an increasing sequence of Hilbertian seminorms $\left(\left\|.\right\|_{k}\right)_{k\in \mathbb{N}}$. For a Fr\'echet space $E$, we  will denote the class of all neighborhoods of zero in $E$ and  the class of all bounded sets in $E$ by $\mathcal{U}\left(E \right)$ and  $\mathcal{B}\left(E \right)$, respectively. If $U$ and $V$ are absolutely convex sets of $E$ and $U$ absorbs $V$, that is, $V\subseteq CU$ for some $C>0$, and $L$ is a subspace of $E$, then we set;
	$$\delta\left(V, U, L\right)=\inf\left\{t>0: V\subseteq tU+L\right\}.$$
	The $n^{th}$ \textit{Kolmogorov diameter} of $V$ with respect to $U$ is defined as;
	$$d_{n}\left(V,U\right)=\inf\left\{\delta\left(V, U, L\right): \dim L\leq n\right\}\hspace{0.2in}n=0,1,2,...$$
	Let $U_{1}\supset U_{2}\supset \cdots \supset U_{p}\supset\cdots$ be a base of neighborhoods of zero of Fr\'echet space E. The diametral dimension of $E$ is defined as
	$$\Delta\left(E \right)= \left\{\left(t_{n}\right)_{n\in \mathbb{N}}: \forall p\in \mathbb{N} \hspace{0.1in} \exists \hspace{0.05in} q>p \hspace{0.1in}\lim_{n\rightarrow\infty} t_{n}d_{n}\left(U_{q},U_{p}\right)=0\right\}.$$
	Demeulenaere et al. \cite{L} showed that the diametral dimension of a nuclear Fr\'echet space can also be represented as
	$$\Delta\left(E\right)=\left\{\left(t_{n}\right)_{n\in \mathbb{N}}: \forall \hspace{0.05in}p\in \mathbb{N} \hspace{0.1in} \exists \hspace{0.05in}q> p \hspace{0.1in}\sup_{n\in \mathbb{N}} \left|t_{n}\right|d_{n}\left(U_{q},U_{p}\right)<+\infty \right\}.$$
	The approximate diametral dimension of a Fr\'echet space $E$ is defined as
	\begin{equation*}
		\delta\left(E \right)= \left\{\left(t_{n}\right)_{n\in \mathbb{N}}: \exists\hspace{0.05in}U\in \mathcal{U}\left(E\right) \hspace{0.05in} \exists \hspace{0.05in}B\in\mathcal{B}\left(E\right) \hspace{0.1in}\lim_{n\rightarrow\infty}\frac{ t_{n}}{d_{n}\left(B,U\right)}=0\right\}.
	\end{equation*}
	It follows from Proposition 6.6.5 of \cite{R} that for a Fr\'echet space $E$ with the base of neighborhoods $U_{1}\supset U_{2}\supset \cdots \supset U_{p}\supset\cdots$, the approximate diametral dimension can be represented as;
	$$\delta\left(E \right)= \left\{\left(t_{n}\right)_{n\in \mathbb{N}}: \exists p\in \mathbb{N} \hspace{0.05in} \forall \hspace{0.05in} q>p \hspace{0.05in}\lim_{n\rightarrow\infty}\frac{ t_{n}}{d_{n}\left(U_{q},U_{p}\right)}=0\right\}.$$
	The following proposition shows how the diametral dimension and the approximate diametral dimension passes into subspaces:
	\begin{prop}\label{inv} Let $E$ be a Fr\'echet space and $F$ be a subspace or a qoutient of $E$. Then,
		~\\ 1. $\Delta\left(E\right)\subseteq \Delta \left(F\right)$.
		~\\ 2. $\delta\left(F\right)\subseteq \delta\left(E\right)$.
		~\\	Hence the diametral dimension and the approximate diametral dimension are linear topological invariants.
	\end{prop}
	\begin{proof} \cite[Proposition 6.6.7 and Proposition 6.6.25]{R}
	\end{proof}
	\par A matrix $\left(a_{k,n}\right)_{k,n\in \mathbb{N}}$ of non-negative numbers is called a \textit{K\"{o}the matrix} if it is satisfies that for each $k\in \mathbb{N}$ there exists an $n\in \mathbb{N}$ with $a_{k,n}>0$ and $a_{k,n}\leq a_{k,n+1}$ for all $k,n\in\mathbb{N}$. For a K\"{o}the matrix $\left(a_{k,n}\right)_{k,n\in \mathbb{N}}$,
	$$\displaystyle K\left(a_{k,n}\right)=\left\{x=\left(x_{n}\right): \left\|x\right\|_{k}:=\sum^{\infty}_{n=1}\left|x_{n}\right|a_{k,n}<+\infty \textnormal{ for all } k\in \mathbb{N}\right\}$$
	is called a \textit{K\"{o}the space}. Every K\"{o}the space is a Fr\'echet space given by the semi-norms in its definition. Nuclearity of a K\"{o}the space was characterized as follows:
	\begin{thm}\label{GP}[Grothendieck-Pietsch] $K\left(a_{kn}\right)$ is nuclear K\"{o}the space if and only if for every $k\in \mathbb{N}$, there exists a $l>k$ so that $\displaystyle \sum^{\infty}_{n=1} \frac{a_{k,n}}{a_{l,n}}<+\infty$.
	\end{thm}
	\begin{proof} \cite[Theorem 28.15]{MV}.
	\end{proof}
	\par Dynin-Mitiagin Theorem \cite[Theorem 28.12]{MV} states that if a nuclear Fr\'echet space E with the sequence of seminorms $\left( \left\|.\right\|_{k}\right) _{k\in \mathbb{N}}$ has a Schauder basis $\left( e_{n}\right)_{n\in \mathbb{N}}$, then it is canonically isomorphic to a nuclear K\"{o}the space defined by the matrix $\left(  \left\|e_{n} \right\|_{k}\right)_{k,n\in \mathbb{N}}$. Therefore, it is important to understand the structure of nuclear K\"{o}the spaces in the theory of nuclear Fr\'echet spaces.
	\par Terzio\u{g}lu gave an estimation for $n^{th}$-Kolmogorov diameters of a K\"{o}the space $K(a_{k,n})$ by using the matrix $(a_{k,n})_{k,n\in \mathbb{N}}$.
	\begin{prop} Let $K(a_{k,n})$ be a K\"{o}the space and fixed $n\in \mathbb{N}$. Assume $J\subset \mathbb{N}$ with $\left|J\right|=n+1$ and $I\subset \mathbb{N}$ with $\left| I\right|\leq n$. Then for every $p$ and $q>p$, 
		$$\displaystyle \inf\left\lbrace \frac{a_{p,i}}{a_{q,i}} :i\in J\right\rbrace \leq d_{n}(U_{q}, U_{p})\leq \sup \left\lbrace \frac{a_{p,i}}{a_{q,i}} :i\notin I  \right\rbrace .$$
	\end{prop}
	\begin{proof} \cite[Proposition 1]{TTtj}.
	\end{proof}
	\begin{defnt}\label{dregularity}
		A K\"{o}the space $K(a_{k,n})$ is called regular if  the inequality 
		$\displaystyle \frac{a_{k+1,n}}{a_{k,n}}\leq \frac{a_{k+1,n+1}}{a_{k,n+1}}$
		is satisfied for all $ k,n\in \mathbb{N}$
	\end{defnt} 
	\begin{remark}\label{nthdiameter}
		In the light of the above proposition, we conclude that for any regular K\"{o}the space $K\left(a_{p,n}\right)$, the $n^{th}$-Kolmogorov diameter is $\displaystyle d_{n}\left(U_{q}, U_{p}\right)=\frac{a_{p, n+1}}{a_{q,n+1}}$. If, on the other hand, $K\left(a_{p,n}\right)$ is not regular, then, one can find Kolmogorov diameters by rewriting the sequence $\displaystyle \left( \frac{a_{p,n}}{ a_{q,n}}\right)_{n\in \mathbb{N}}$ with terms in a descending order so that the $n^{th}$-Kolmogorov diameter of $K\left(a_{p,n}\right)$ is nothing but the $n+1-th $ term of this descending sequence.
	\end{remark}
	\par Power series spaces are the most important family of Köthe spaces and they have a significant role in this work, for a comprehensive survey see \cite{TTV}. Let $\alpha=\left(\alpha_{n}\right)_{n\in \mathbb{N}}$ be a non-negative increasing sequence with $\displaystyle \lim_{n\rightarrow \infty} \alpha_{n}=+\infty$. A power series space of finite type is defined by
	$$\displaystyle \Lambda_{1}\left(\alpha\right):=\left\{x=\left(x_{n}\right)_{n\in \mathbb{N}}: \left\|x\right\|_{k}=\sum^{\infty}_{n=1}\left|x_{n}\right|e^{\displaystyle \scalebox{0.95}{$-{\displaystyle \scalebox{1.1}{$\frac{1}{k}$}} \alpha_{n}$}}<+\infty \textnormal{ for all } k\in \mathbb{N}\right\}$$
	and a power series space of infinite type is defined by
	$$\displaystyle \Lambda_{\infty}\left(\alpha\right):=\left\{x=\left(x_{n}\right)_{n\in \mathbb{N}}: \left\|x\right\|_{k}=\sum^{\infty}_{n=1}\left|x_{n}\right|e^{\displaystyle \scalebox{0.9}{$\hspace{0.025in}k\alpha_{n}$}}<+\infty \textnormal{ for all } k\in \mathbb{N}\right\}.$$
	~\\$\vspace{-.35in}$~\\ The nuclearity of a power series space of finite type $\Lambda_{1}\left(\alpha\right)$ and of infinite type $\Lambda_{\infty}\left(\alpha\right)$ are equivalent to the conditions $\displaystyle \lim_{n\rightarrow \infty} \frac{\ln(n)}{\alpha_{n}}=0$ and $\displaystyle \sup_{n\in \mathbb{N}} \frac{\ln(n)}{\alpha_{n}}<+\infty$, respectively. 
	\begin{defnt}An exponent sequence $\alpha$ is called {finitely nuclear} if $\Lambda_{1}(\alpha)$ is nuclear.
	\end{defnt}
	~\\$\vspace{-.4in}$~\\Diametral dimension and approximate diametral dimension of power series spaces are $\displaystyle\Delta\left(\Lambda_{1}\left(\alpha\right)\right)=\Lambda_{1}\left(\alpha\right)$, $\hspace{0.05in}\Delta\left(\Lambda_{\infty}\left(\alpha\right)\right)=\Lambda_{\infty}\left(\alpha\right)^{\prime}$, $\hspace{0.05in}\delta\left(\Lambda_{1}\left(\alpha\right)\right)=\Lambda_{1}\left(\alpha\right)^{\prime}$ and $\hspace{0.05in}\delta\left(\Lambda_{\infty}\left(\alpha\right)\right)=\Lambda_{\infty}\left(\alpha\right)$ for details see \cite{BRP} and \cite{M2}.
	\par An exponent sequence $\alpha$ is called
	\begin{center}
		\begin{itemize}
			\item[]{}\textit{$\hspace{1in}$stable} $\hspace{.8in}$if $\hspace{0.25in}\displaystyle \sup_{n\in \mathbb{N}}\hspace{0.025in}\frac{\alpha_{\hspace{0.01in}2n}}{ \alpha_{n}}<+\infty$, ~\\~\\
			\item[]{}\textit{$\hspace{1in}$weakly-stable} $\hspace{0.25in}$ if $ \hspace{0.25in}\displaystyle \sup_{n\in \mathbb{N}}\hspace{0.025in}\frac{\alpha_{\hspace{0.01in}n+1}}{\alpha_{n}}<+\infty$, $\hspace{0.15in}$~\\~\\ 
			\item[]{}\textit{$\hspace{1in}$unstable}  $\hspace{0.65in}$if  $\hspace{0.25in}\displaystyle \lim_{n\rightarrow\infty}\hspace{0.025in}\frac{\alpha_{\hspace{0.01in}n+1}}{\alpha_{n}}=+\infty$.
		\end{itemize}
	\end{center}
	It follows that $\alpha$ is stable, respectively weakly-stable, if and only if $E\cong E\times E$, respectively, $E\cong E\times \mathbb{K}$ where $E=\Lambda_{r}(\alpha)$ for $r=1$ or $r=\infty$, for proofs see \cite{Du}. 
	\par Subspaces and quotient spaces of a nuclear
	stable power series space are characterized by Vogt and Wagner (\cite{VDN},\cite{VW}, \cite{VW2}) in terms of diametral dimension and DN-$\Omega$ type linear-topological invariants. The topological invariants \underline{DN} and $\Omega$ are enjoyed by many natural nuclear Fr\'echet
	spaces appearing in analysis and these invariants play an important role in this study.
	~\\$\vspace{-.45in}$~\\
	\begin{defnt}
		A Fr\'echet space $(E, \left\|.\right\|_{k})_{k\in \mathbb{N}}$ is said to have the property:
		~\\$\vspace{-.15in}$~\\
		{\textbf{(\underline{DN})}}$\hspace{0.1in} \exists \hspace{0.05in}k $ $\hspace{0.05in} \forall \hspace{0.05in} j $ $\hspace{0.05in} \exists \hspace{0.05in}l,\hspace{0.05in} C>0, \hspace{0.05in} 0<\lambda <1$  
		$$\hspace{1.65in} \displaystyle \left\| x\right\|_{j}\leq C\hspace{0.025in} \left\|x \right\|^{\lambda}_{k} \hspace{0.025in} \left\|x \right\|^{1-\lambda}_{l}  \hspace{1.25in} \forall\hspace{0.05in} x\in E $$
		~\\$\vspace{-.325in}$~\\
		{\ka } $\hspace{0.2in} \forall \hspace{0.05in}p$ $\hspace{0.05in} \exists \hspace{0.05in} q $ $\hspace{0.05in} \forall \hspace{0.05in}k$ $\hspace{0.05in} \exists \hspace{0.05in} C>0, \hspace{0.05in} 0<\tau <1$  
		$$\hspace{1.65in} \displaystyle {\left\|y\right\|_{q}^{*}}\leq C{\left\|y\right\|_{p}^{*}}^{1-\theta}{\left\|y\right\|_{k}^{*}}^{\theta}\hspace{1.25in} \forall \hspace{0.05in} y\in E^{'}$$ ~\\$\vspace{-.325in}$~\\ where $\displaystyle \left\|y\right\|^{*}_{k}:=\sup \left\{ \left|y\left(x\right)\right|: \left\|x\right\|_{k} \leq 1\right\}\in \mathbb{R}\cup \left\{+\infty\right\}$	is the gauge functional of the polar $U^{\circ}_{k}$ for $U_{k}= \left\{x\in E: \left\|x\right\|_{k}\leq 1\right\}$. \end{defnt}
	~\\$\vspace{-.3in}$~\\
	In \cite{VSQ}, D. Vogt characterized $\Omega$ for K\"{o}the spaces in terms of K\"{o}the matrix as follows:
	\begin{prop}\label{pomega} A K\"{o}the space $K\left(a_{k,n}\right)$ has the property $\Omega$ if and only if the condition
		$$\displaystyle \forall\hspace{0.05in} p\hspace{0.1in} \exists\hspace{0.05in} q \hspace{0.1in}\forall\hspace{0.05in}k \hspace{0.1in}\exists \hspace{0.05in}j>0, \hspace{0.05in}C>0 \hspace{0.5in} { \left(a_{p,n}\right)^{j}a_{k,n}\leq C\left(a_{q,n}\right)^{j+1}} \hspace{0.5in} \forall n\in\mathbb{N}$$
		is satisfied.
	\end{prop} 
	\begin{proof} \cite[Proposition 5.3]{VSQ}.
	\end{proof}
	By using the technique in \cite[5. 1 Proposition]{VSQ}, one can easily obtain the following:
	\begin{prop}\label{pdn} A K\"{o}the space $K\left(a_{k,n}\right)$ has the property $\underline{DN}$ if and only if the condition
		$$\displaystyle \exists\hspace{0.05in} p_{0}\hspace{0.1in} \forall\hspace{0.05in} p \hspace{0.1in}\exists \hspace{0.05in}q \hspace{0.1in}\exists \hspace{0.05in}0<\lambda<1, \hspace{0.05in}C>0 \hspace{0.35in} a_{p,n}\leq C  \left(a_{p_{0},n}\right)^{\lambda}\left(a_{q,n}\right)^{1-\lambda} \hspace{0.2in} \forall n\in\mathbb{N}$$
		is satisfied.
	\end{prop} \par
	Now we give the important result which gives a relation between the diametral dimension/approximate diametral dimension of a nuclear Fr\'echet spaces with the properties $\underline{DN}$, $\Omega$ and that of a power series spaces $\Lambda_{1}\left(\varepsilon\right)$ and $\Lambda_{\infty}\left(\varepsilon\right)$ for some special exponent sequence $\varepsilon$.
	\begin{prop} Let $E$ be a nuclear Fr\'echet space with the properties $\underline{DN}$ and $\Omega$. There exists an exponent sequence (unique up to equivalence) $\varepsilon=\left(\varepsilon_{n}\right)$ satisfying:
		\begin{equation}\Delta\left(\Lambda_{1}\left(\varepsilon\right)\right)\subseteq\Delta\left(E\right)\subseteq\Delta\left(\Lambda_{\infty}\left(\varepsilon\right)\right). 
		\end{equation}
		Furthermore, $\Lambda_{1}\left(\alpha\right)\subseteq\Delta\left(E\right)$ implies $\Lambda_{1}\left(\alpha\right)\subseteq \Lambda_{1}\left(\varepsilon\right)$ and $\Delta\left(E\right)\subseteq \Lambda^{\prime}_{\infty}\left(\alpha\right)$ implies $\Lambda^{\prime}_{\infty}\left(\varepsilon\right)\subseteq\Lambda^{'}_{\infty}\left(\alpha\right)$.
	\end{prop}
	\begin{proof} \cite[Proposition 1.1]{AKT2}.
	\end{proof}
	\begin{defnt}\label{aes} Let $E$ be a nuclear Fr\'echet space with the properties $\underline{DN}$ and $\Omega$. The sequence $\varepsilon$ (unique up to equivalence) in the above proposition is called the \textbf{associated exponent sequence} of $E$ in \cite{AKT2}. 
	\end{defnt}
	We note that $\Lambda_{\infty} (\varepsilon)$ is always nuclear provided $E$ is nuclear, but it may happen that $\Lambda_{1} (\varepsilon)$ is not nuclear. For example, if we take the space of rapidly decreasing sequence $s=\Lambda_{\infty}(\ln(n))$, the associated exponent sequence of $s$ is $(\ln(n))_{n\in \mathbb{N}}$ and $\Lambda_{1}(\ln(n))$ is not nuclear.
	\par In the proof of the above proposition, Aytuna et al. showed that there exists an exponent sequence (unique up to equivalence) $\left(\varepsilon_{n}\right)$ such that for each $p\in \mathbb{N}$ and $q>p$, there exist $C_{1}, C_{2}>0$ and $a_{1}, a_{2}>0$ satisfying
	$$C_{1}\hspace{0.025in}e^{\displaystyle {-\scalebox{1}{$\displaystyle a_{\scalebox{0.65}{$1$}}\hspace{0.025in}\varepsilon_{n}$}}}\leq d_{n}\left(U_{q},U_{p}\right)\leq C_{2}\hspace{0.025in}e^{\displaystyle {-\scalebox{1}{$\displaystyle a_{\scalebox{0.65}{$2$}}\hspace{0.025in}\varepsilon_{n}$}}}$$
	for all $n\in \mathbb{N}$. From this inequality, one can easily obtain
	\begin{center} 
		$ \delta\left(\Lambda_{\infty}\left(\varepsilon\right)\right)\subseteq\delta\left(E\right)\subseteq\delta\left(\Lambda_{1}\left(\varepsilon\right)\right).
		$\end{center}
	\par For a nuclear Fr\'echet space $E$ with the properties $\underline{DN}$ and $\Omega$ and the associated exponent sequence $\varepsilon$, concidence of the diametral dimension of $E$ with that of power series spaces defined by $\varepsilon$ form two extreme cases. The extreme case $\Delta(E)=\Delta\left( \Lambda_{\infty}(\varepsilon)\right)$ gives an information about a (complemented) subspace  of a nuclear Fr\'echet space $E$ with the properties $\underline{DN}$ and $\Omega$ and stable associated exponent sequence $\varepsilon$. In \cite{AKT2}, Aytuna et al. proved that a nuclear Fr\'echet space $E$ with the properties $\underline{DN}$ and $\Omega$ contains a complemented copy of $\Lambda_{\infty}(\varepsilon)$ provided that $\Delta(E)=\Delta(\Lambda_{\infty}(\varepsilon))$ and $\varepsilon$ is stable.
	\begin{thm} Let $E$ be a nuclear Fr\'echet space with the properties $\underline{DN}$ and $\Omega$ and stable associated exponent sequence $\varepsilon$. If $\Delta\left(E\right)=\Delta\left(\Lambda_{\infty}\left(\varepsilon\right)\right)$, then $E$ has complemented subspace which is isomorphic to $\Lambda_{\infty}\left(\varepsilon\right)$.
	\end{thm}
	\begin{proof} \cite[Theorem 1.2]{AKT2}.
	\end{proof}
	\par On the other hand, there is no information for the other extreme $\Delta(E)=\Delta\left( \Lambda_{1}(\varepsilon)\right)$. This leads to ask the Question \ref{qi1} in Introduction. We need the following proposition characterizing coincidence of $\delta(E)$ with $\delta(\Lambda_{1}(\varepsilon))$ given by A. Aytuna in \cite{A}: 
	\begin{prop}\label{AA} Let $E$ be a nuclear Fr\'echet space $E$ with the properties $\underline{DN}$, $\Omega$ and associated exponent sequence $\varepsilon$. Then
		$$\delta\left(E\right)= \delta\left(\Lambda_{1}\left(\varepsilon\right)\right) \hspace{0.25in}\Leftrightarrow \hspace{0.25in} \inf_{p}\sup_{q\geq p} \limsup_{n\in \mathbb{N}}\frac{\varepsilon_{n}\left(p,q\right)}{\varepsilon_{n}}=0$$
		where $\varepsilon_{n}\left(p,q\right)=-\log{d_{n}\left(U_{q},U_{p}\right)}$. 
	\end{prop}
	\begin{proof} \cite[Corollary 1.10]{A}
	\end{proof}
	\section{ $\bf \mathcal{K}_{\alpha}$ Spaces}
	\par In this section, we will construct a family of nuclear K\"{o}the spaces with the properties $\underline{DN}$ and $\Omega$ and parameterized by a finitely nuclear sequence $\alpha$ and show that a subfamily of these K\"{o}the spaces satisfied that  $\Delta\left(K\left(a_{k,n}\right)\right)=\Delta\left( \Lambda_{1}\left(\varepsilon\right)\right) $ and $\delta\left(K\left(a_{k,n}\right)\right)\neq \delta\left( \Lambda_{1}\left(\varepsilon\right)\right)$ for its associated exponent sequence $\varepsilon$. This shows that Question \textbf{\ref{qi2}}  has a negative answer. 
	\par We proceed as follows: First, we divide natural numbers $\mathbb{N}$ into infinite disjoint union of infinite subsets. For this purpose, we order the elements of $\mathbb{N}^{2}$ by matching them with the elements of $\mathbb{N}$ such that any element  $\left(x,y\right)\in \mathbb{N}^{2}$ corresponds to the element $\displaystyle \frac{\left(x+1\right)\left(x+2\right)}{2} +y(x+1)+ \frac{y\left(y-1\right)}{2}\in \mathbb{N}$. One can visualize this ordering as shown in the following diagram:
	\begin{center}
		\scalebox{0.75}
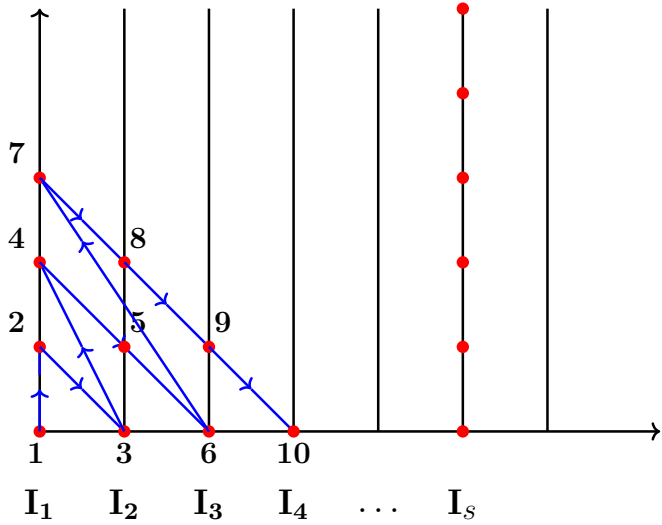
{
			\begin{tikzpicture}[point/.style={fill,circle,inner sep=2.2pt,}]
				\draw [<->,line width=0.45mm ] (0,7.5) -- (0,0) -- (11,0);
				\draw[,line width=0.45mm,](1.5,0)--node[above=3mm,right]{}(1.5,7.5); 
				
				\draw[-,line width=0.45mm, ](3,0)--node[above=3mm,right]{}(3,7.5); 
				
				\draw[-,line width=0.45mm, ](4.5,0)--node[above=3mm,right]{}(4.5,7.5); 
				\draw[-,line width=0.45mm, ](6,0)--node[above=3mm,right]{}(6,7.5); 
				\draw[-,line width=0.45mm, ](7.5,0)--node[above=3mm,right]{}(7.5,7.5); 
				\draw[-,line width=0.45mm, ](9,0)--node[above=3mm,right]{}(9,7.5); 
				
				
				
				

				\node[,label={${\hspace{-0.25in}\scalebox{1.3}{\bf{1}}}$}] at (0.25, -0.85){};
				\node[point,red,label={}] at (0,0){};
				\draw[->,line width=0.45mm,blue](0,0)--node[above=3mm,right]{}(0,0.75); 
				\draw[-, line width=0.45mm,blue](0,1)--node[above=3mm,right]{}(0,1.5); 
				\node[point,red,label={$\scalebox{1.3}{\bf{\hspace{-0.25in}2}}$}] at (0,1.5){};
				\draw[->,line width=0.45mm,blue](0,1.5)--node[above=3mm,right]{}(0.75,0.75); 
				\draw[-, line width=0.45mm,blue](0.75,0.75)--node[above=3mm,right]{}(1.5,0); 
				\node[point,red, label={}] at (1.5,0){};
				\node[label={\scalebox{1.3}{$\bf{3}$}}] at (1.5,-0.85){};
				\draw[->,line width=0.45mm,blue](1.5,0)--node[above=3mm,right]{}(0.75,1.5); 
				\draw[-, line width=0.45mm,blue](0.75,1.5)--node[above=3mm,right]{}(0,3); 
				\node[point,red,label={$\scalebox{1.3}{\bf{\hspace{-0.25in}4}}$}] at (0,3){};
				\draw[->,line width=0.45mm,blue](0,3)--node[above=3mm,right]{}(1.5,1.5); 
				\draw[-, line width=0.45mm,blue](1.5,1.5)--node[above=3mm,right]{}(3,0); 
				\node[point,red,label={$\scalebox{1.3}{\bf{\hspace{0.15in}5}}$}] at (1.5,1.5){};
				\node[point,red, label={}] at (3,0){};
				\node[label={$\scalebox{1.3}{\bf{6}}$}] at (3,-0.85){};
				\draw[->,line width=0.45mm,blue](3,0)--node[above=3mm,right]{}(0.75,3.375); 
				\draw[-, line width=0.45mm,blue](0.75,3.375)--node[above=3mm,right]{}(0,4.5); 
				\node[point,red,label={$\scalebox{1.3}{\bf{\hspace{-0.25in}7}}$}] at (0,4.5){};
				\draw[->,line width=0.45mm,blue](0,4.5)--node[above=3mm,right]{}(0.75,3.75); 
				\draw[-, line width=0.45mm,blue](0.75,3.75)--node[above=3mm,right]{}(1.5,3); 
				\node[point,red,label={$\scalebox{1.3}{\bf{\hspace{0.15in}8}}$}] at (1.5,3){};
				\draw[->,line width=0.45mm,blue](1.5,3)--node[above=3mm,right]{}(2.25,2.25); 
				\draw[-, line width=0.45mm,blue](2.25,2.25)--node[above=3mm,right]{}(3,1.5); 
				\node[point,red,label={$\scalebox{1.3}{\bf{\hspace{0.15in}9}}$}] at (3,1.5){};
				\draw[->,line width=0.45mm,blue](3,1.5)--node[above=3mm,right]{}(3.75,0.75); 
				\draw[-, line width=0.45mm,blue](3.75,0.75)--node[above=3mm,right]{}(4.5,0); 
				\node[point,red,label={}] at (4.5,0){};
				\node[label={$\scalebox{1.3}{\bf{10}}$}] at (4.5,-0.85){};
				
				\node[label={\scalebox{1.5}{$\bf \displaystyle {{I_{{1}}}}$}}] at (0,-1.85){};
				\node[label={\scalebox{1.5}{$ \bf\displaystyle {{I_{{2}}}}$}}] at (1.5,-1.85){};
				\node[label={\scalebox{1.5}{$\bf \displaystyle {{I_{{3}}}}$}}] at (3,-1.85){};
				\node[label={\scalebox{1.5}{$\bf \displaystyle {{I_{{4}}}}$}}] at (4.5,-1.85){};
				\node[label={\scalebox{1.5}{$\bf\displaystyle I_{\displaystyle \scalebox{0.85}{$s$}}$}}] at (7.5,-1.85){};
				\node[label={\scalebox{1.5}{$\bf \displaystyle {\mathbf{\cdots}}$}}] at (6,-1.85){};
				\node[point,red,label={}] at (7.5,0){};
				\node[point,red,label={}] at (7.5,1.5){};
				\node[point,red,label={}] at (7.5,3){};
				\node[point,red,label={}] at (7.5,4.5){};
				\node[point,red,label={}] at (7.5,6){};
				\node[point,red,label={}] at (7.5,7.5){};
		\end{tikzpicture}}
	\end{center}
	As shown in the above diagram, each vertical line $I_{\displaystyle \scalebox{0.85}{$s$}}$ has infinitely many elements and $\mathbb{N}$ can be expressed as an infinite disjoint union of $I_{\displaystyle \scalebox{0.85}{$s$}}$, that is, $\displaystyle N=\bigcup_{s\in \mathbb{N}}I_{\displaystyle \scalebox{0.85}{$s$}}$. 
	\begin{defnt}
		Let $\displaystyle\alpha= \left(\alpha_{n}\right)_{n\in \mathbb{N}}$ be a strictly increasing, positive, finitely nuclear sequence. We define a matrix $\displaystyle \left( a_{k,n}\right) _{k,n\in\mathbb{N}}$ by setting:
		\begin{equation}\label{es}
			\displaystyle a_{k,n}= 
			\begin{cases}
				e^{\displaystyle -\frac{1}{ k}\hspace{0.025in}\alpha_{n}},& \text{if}\hspace{0.1in} k\leq s \\ \\
				e^{\displaystyle \left(-\frac{1}{k}+1\right)\alpha_{n}},& \text{if} \hspace{0.1in}  k\geq s+1.
			\end{cases}
		\end{equation}
		where $n\in I_{\displaystyle \scalebox{0.85}{$s$}}$,  $s\in \mathbb{N}$.
	\end{defnt}
	Infact, $\displaystyle \left(a_{k,n}\right)_{k,n\in \mathbb{N}}$ is a K\"{o}the matrix, since for every  $n,k\in \mathbb{N}$,    $\hspace{0.05in}\displaystyle 0< a_{k,n}\leq a_{k+1,n}$. We denote the K\"{o}the space generated by a matrix $\left(a_{k,n}\right)_{k,n\in \mathbb{N}}$ as in \ref{es} by $\bf \mathcal{K}_{\alpha}$. We say that the space $\bf \mathcal{K}_{\alpha}$ is parameterized by the sequence $\alpha$. We denote the family of all Köthe space $\bf \mathcal{K}_{\alpha}$ by $\bf \mathcal{K}$.
	Now, we show that each element of the family $\bf \mathcal{K}$ is nuclear and satisfies the properties $\underline{DN}$ and $\Omega$:
	\begin{lemma} Let $\bf \mathcal{K}_{\alpha}$ be an element of the family $\bf \mathcal{K}$ parametrized by $\alpha= \left(\alpha_{n}\right)_{n\in \mathbb{N}}$. Then, $\bf \mathcal{K}_{\alpha}$ is nuclear and has the properties  $\underline{DN}$ and $\Omega$.
	\end{lemma}
	\begin{proof} For the nuclearity of $\bf \mathcal{K}_{\alpha}$, we show that the series $\displaystyle \mathlarger{\sum}^{\infty}_{n=1} \hspace{0.025in}{a_{k,n}\over a_{k+1,n}}$ is convergent for each $k\in \mathbb{N}$. Since $\displaystyle \frac{a_{k,n}}{a_{k+1,n}}\leq e^{\mathlarger{\left(-\frac{1}{k}+\frac{1}{k+1}\right)\alpha_{n}}}$ for every $k,n\in \mathbb{N}$ and $\Lambda_{1}(\alpha)$ is nuclear, then the series  $\displaystyle \mathlarger{\sum}^{\infty}_{n=1}\hspace{0.025in} \frac{a_{k,n}}{a_{k+1,n}}$ is convergent. By Theorem \ref{GP}, $\bf \mathcal{K}_{\alpha}$ is nuclear, as asserted.
		\par We now prove that $\bf \mathcal{K}_{\alpha}$ has the $\underline{DN}$ property by using Proposition \ref{pdn}. We will show that for all $p\in \mathbb{N}$ there exists a $0<\lambda<1$ such that the inequality \begin{equation}\label{e11}
			\displaystyle {a_{p,n}\leq  \left(a_{1,n}\right)^{\lambda}\left(a_{p+1,n}\right)^{1-\lambda}}
		\end{equation} 
		is satisfied for all $n\in \mathbb{N}$. Let $p,n\in \mathbb{N}$ and assume $n\in I_{\displaystyle \scalebox{0.85}{$s$}}$, $s\in \mathbb{N}$. There are two cases for $p$ and $s$: $p\leq s$ or $p>s$. First we assume that $p\leq s$: In this case, $a_{1,n}=e^{\displaystyle {-\alpha_{n}}}$, $\displaystyle a_{p,n}=e^{\mathlarger{-{1\over p}}\hspace{0.025in}{\displaystyle \alpha_{n}}}$ and $ a_{p+1,n}\geq e^{\mathlarger{-{1\over p+1}}\hspace{0.025in}{\displaystyle \alpha_{n}}}.$ 
		Then, the inequality \ref{e11} is satisfied for any $\displaystyle  \lambda<\frac{ \frac{1}{ p}-\frac{1}{p+1} }{1-\frac{1}{p+1} }.$ Second we assume that $s<p$: In this case,
		$a_{1,n}=e^{\displaystyle {-\alpha_{n}}}$, $a_{p,n}=e^{\left({\mathlarger{-\frac{1}{p}}}+1\right){\displaystyle \alpha_{n}}}$ and $a_{p+1,n}=e^{\left({\mathlarger{-\frac{1}{p+1}}}+1\right){\displaystyle \alpha_{n}}}.$
		But then the inequality \ref{e11} is satisfied for any $\displaystyle \lambda<\frac{\frac{1}{p}-\frac{1}{p+1}}{2-\frac{1}{p+1}}.$
		Hence, if we choose a $\lambda>0$ satisfying  $$\displaystyle \lambda < \min\left\lbrace  {\frac{\frac{1}{p}-\frac{1}{p+1}}{1-\frac{1}{p+1}}},{\frac{\frac{1}{p}-\frac{1}{p+1}}{2-\frac{1}{p+1}}} \right\rbrace = {\frac{\frac{1}{p}-\frac{1}{p+1}}{2-\frac{1}{p+1}}}$$ 
		then inequality \ref{e11} holds in general and so $\bf \mathcal{K}_{\alpha}$ has the property \underline{DN}, as claimed. 
		\par We now prove that $\bf \mathcal{K}_{\alpha}$ has $\Omega$ by using Proposition \ref{pomega}. We will show that for all $p\in \mathbb{N}$ and $k>p$ there exists a $j>0$ such the inequality 
		\begin{equation}\label{eq2}{ \left(a_{p,n}\right)^{j}a_{k,n}\leq \left(a_{p+1,n}\right)^{j+1}}\end{equation}
		is satisfied for all $n\in \mathbb{N}$. Let $p,n\in \mathbb{N}$ and assume $n\in I_{\displaystyle \scalebox{0.85}{$s$}}$, $s\in \mathbb{N}$. There are two case for $p$ and $s$: $p\leq s$ or $p>s$. First we assume that $p\leq s$: In this case, $a_{p,n}=e^{\mathlarger{-\frac{1}{ p}}\hspace{0.01in}{\displaystyle\alpha_{n}}}$, $a_{p+1,n}\geq e^{\mathlarger{-\frac{1}{ p+1}}\hspace{0.01in}{\displaystyle \alpha_{n}}}$ and $ a_{k,n}\leq e^{\left({\mathlarger{-\frac{1}{ k}}+1}\right){\displaystyle \alpha_{n}}}$
		for all $k\geq p$. Then, the inequality \ref{eq2} is satisfied for any $\displaystyle j\geq {\frac{\frac{1}{p+1}-\frac{1}{k}+1} {\frac{1}{p}-\frac{1}{p+1}}}.$ Second we assume that $s<p$: In this case, $a_{p,n}=e^{\left({\mathlarger{-\frac{1}{ p}}}+1\right){\displaystyle \alpha_{n}}}$, $ a_{p+1,n}=e^{\left({\mathlarger{-\frac{1}{ p+1}}+1}\right){\displaystyle \alpha_{n}}}$ and $ a_{k,n}=e^{\left({\mathlarger{-\frac{1}{ k}}+1}\right){\displaystyle \alpha_{n}}}$
		for all $k\geq p$. Therefore, the inequality \ref{eq2} is satisfied for any $\displaystyle j\geq {\frac{\frac{1}{p+1}-\frac{1}{k}} {\frac{1}{p}-\frac{1}{p+1}}}.$ Now, we choose a $j>0$ satisfying 
		$$\displaystyle j\geq \max\left( {\frac{\frac{1}{ p+1}-\frac{1}{k}+1}{\frac{1}{ p}-\frac{1}{p+1}}}, {\frac{\frac{1}{p+1}-\frac{1}{k}}{\frac{1}{ p}-\frac{1}{p+1}}}\right)={\frac{\frac{1}{ p+1}-\frac{1}{k}+1}{\frac{1}{ p}-\frac{1}{p+1}}}$$
		and so that the inequality \ref{eq2} is satisfied for all $n\in \mathbb{N}$. Hence $\bf \mathcal{K}_{\alpha}$  has the property ${\Omega}$, as claimed. 
	\end{proof}
	\begin{remark}\label{rd2} It is worth noting that any element $\bf \mathcal{K}_{\alpha}$ of the family $\bf \mathcal{K}$ does not have the property $(d_{2})$,
		$$(d_{2}):\hspace{0.9in}\displaystyle \forall k\hspace{0.15in} \exists j \hspace{0.15in} \forall l \hspace{0.9in} \sup_{n}\hspace{0.025in}{a_{kn}\hspace{0.025in}a_{ln}\over (a_{jn})^{2}}<+\infty.\hspace{1.2in}$$
		Since for all $j\in \mathbb{N}$, $n\in I_{j}$, $\displaystyle a_{1,n}=e^{\displaystyle{-\alpha_{n}}}$, $a_{jn}=e^{{\mathlarger{-\frac{1}{j}} {\displaystyle\alpha_{n}}}}$, $ a_{j+1,n}=e^{\left(-{\mathlarger{\frac{1}{ j+1}}}+1\right){\displaystyle \alpha_{n}}}$,
		
		$$\displaystyle \frac{a_{1,n}\hspace{0.025in}a_{j+1,n}}{ (a_{jn})^{2}}= e^{{\mathlarger{\frac{j+2}{ j(j+1)}} \displaystyle\alpha_{n}}} \hspace{0.35in} \textnormal{and} \hspace{0.35in} \sup_{n\in I_{j}}\frac{a_{1,n}\hspace{0.025in}a_{j+1,n}}{ (a_{jn})^{2}}=\sup_{n\in \mathbb{N}}\frac{a_{1,n}\hspace{0.025in}a_{j+1,n}}{(a_{jn})^{2}}=+\infty
		$$
		then $\bf \mathcal{K}_{\alpha}$ does not have the property $(d_{2})$. So the family $\bf \mathcal{K}$ does not contain a power series space of finite type.
	\end{remark}
	\subsection{Kolmogorov diameters of an element $\bf \mathcal{K}_{\alpha}$ of the family $\bf \mathcal{K}$}
	~\\$\vspace{-.2in}$  ~\\ \par In this subsection, we calculate Kolmogorov diameters of an element $\bf \mathcal{K}_{\alpha}$ \\of the family $\bf \mathcal{K}$. In order to determine $n^{th}$-Kolmogorov diameter of a K\"{o}the space $\bf \mathcal{K}_{\alpha}$, we will rewrite the sequence $\displaystyle \left(\frac{a_{p,n}}{a_{q,n}}\right)_{n\in \mathbb{N}}$ in descending order. We know from Remark \ref{nthdiameter} that the $n^{th}$-Kolmogorov diameter of the space $\bf \mathcal{K}_{\alpha}$ is the ${n+1}^{th}$-term of this descending sequence.
	\par Let $\bf \mathcal{K}_{\alpha}$ be an element of the family $\bf \mathcal{K}$ parameterized by an exponent sequence $\alpha$. Let us take a $p$, a $q> p$ and an $n\in I_{\displaystyle \scalebox{0.85}{$s$}}$, $s\in \mathbb{N}.$ Then, we can write
	\begin{equation*}\label{pq2}
		\displaystyle {a_{p,n}\over a_{q,n}}= 
		\begin{cases}
			e^{\displaystyle c_{pq}\hspace{0.025in}\alpha_{n}},& \text{}\hspace{0.1in} s\geq q\hspace{0.05in}  \text{or}\hspace{0.05in}  s<p \\ \\
			e^{\displaystyle \left(c_{pq}-1\right)\alpha_{n}},& \text{}\hspace{0.1in} p\leq s< q \\
		\end{cases}
	\end{equation*}
	where $\displaystyle c_{pq}$ is the negative number $\displaystyle  -{1\over p}+{1\over q}$. We define the set $ \displaystyle I=\bigcup_{p\leq s<q} I_{\displaystyle \scalebox{0.85}{$s$}}$ with the elements $(n_{\displaystyle \scalebox{0.85}{$i$}})_{i \in \mathbb{N}}$ ordered increasingly, namely, $n_{i}\leq n_{i+1}$  for all $i\in \mathbb{N}$. We also denote the index of the element of $I_{p}$ on the line with the equation $x+y=q+k-2$ by $\displaystyle \boldsymbol{ \scalebox{1.05}{$s_{\displaystyle\scalebox{0.8} {$k$}}$}}$ for each $k=0,1,2,...$, as seen from the following diagram. Since every a line with the equation  $x+y=q+k-2$ has $q-p$ elements of $I$, then $\displaystyle \scalebox{1.05}{$s_{\displaystyle\scalebox{0.75} {$(k+1)$}}$}- \scalebox{1.05}{$s_{\displaystyle\scalebox{0.8} {$k$}}$}=q-p$ for every $k=0,1,2,....$
	\begin{center}\label{I}
		\scalebox{0.75}
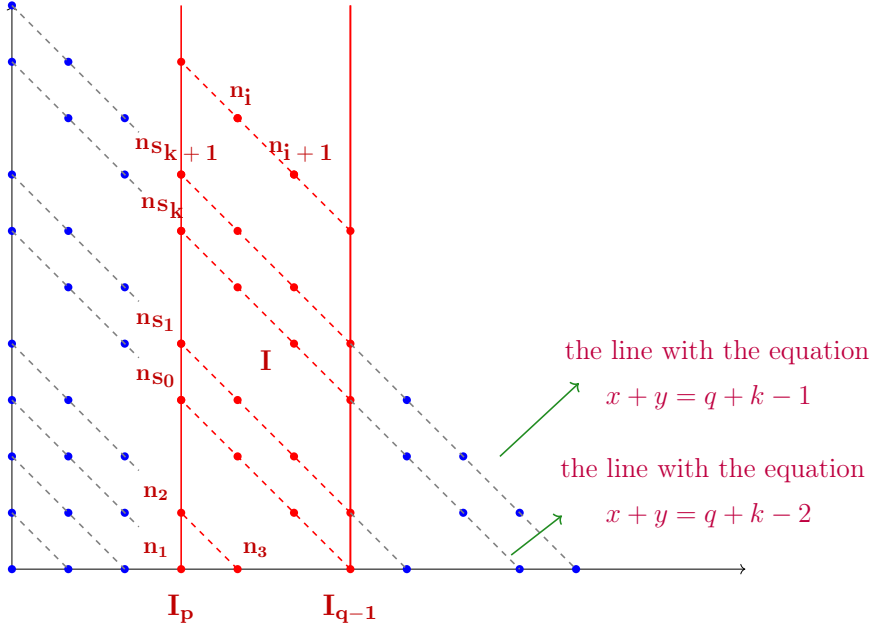
{
			\begin{tikzpicture}[point/.style={fill,circle,inner sep=1.45pt,black}]
				\draw [<->] (1,10) -- (1,0) -- (14,0);
				\draw[-,thick,red](4,0)--node[above=3mm,right]{}(4,10); 
				\draw[-,thick,red](7,0)--node[above=3mm,right]{}(7,10); 
				\draw[-,thick,red](7,0)--node[above=3mm,right]{}(7,10); 

				\draw[dashed,thick,red](5,0)--node[above=3mm,right]{}(4,1); 
				\draw[dashed,thick,red](7,0)--node[above=3mm,right]{}(4,3); 
				
				\node[point,red,label={\color{Mycolor2} $\bf {\hspace{-0.35in} n_{\displaystyle \hspace{0.01in}\scalebox{1.075}{$\bf s_{\displaystyle\scalebox{0.75} {$\bf 0$}}$}}}$}] at (4,3){};
				\node[point,red,label={\color{Mycolor2} $\bf {\hspace{-0.35in} n_{\displaystyle \hspace{0.01in}\scalebox{1.075}{$\bf s_{\displaystyle\scalebox{0.75} {$\bf 1$}}$}}}$}] at (4,4){};

				\draw[dashed,thick,red](7,1)--node[above=3mm,right]{}(4,4); 
				
				\draw[dashed,thick,red](7,3)--node[above=3mm,right]{}(4,6); 
				
				\draw[dashed,thick,red](7,4)--node[above=3mm,right]{}(4,7); 
				
				
				\draw[dashed,thick,red](7,6)--node[above=3mm,right]{}(4,9); 
				
				
				
				
				\node[point,red] at (6,7){};
				\node[point,red] at (7,6){};
				
				\node[point,red] at (5,8){};
				
				\node[point,red] at (4,7){};
				\node[point,red,label={\color{Mycolor2} $\bf {\hspace{-0.05in} n_{\displaystyle \hspace{0.01in}\scalebox{1.075}{$\bf s_{\displaystyle\hspace{0.01in}\scalebox{0.9} {$\bf k+1$}}$}}}$}] at (4,7){};
				
				\node[point,red] at (4,9){};
				
				\node[point,red,label={$\bf \color{Mycolor2}{\hspace{0.05in} \displaystyle n_{\displaystyle i}}$}] at (5,8){};
				\node[point,red,label={$\bf \color{Mycolor2}{\hspace{0.1in} \displaystyle n_{\displaystyle i+1}}$}] at (6,7){};
				
				\node[point,red,label={$\bf \color{Mycolor2}{\hspace{-0.35in}n_{1}}$}] at (4,0){};	
				\node[point,red, label={$\bf \color{Mycolor2}{\hspace{-0.35in}n_{2}}$}] at (4,1){};
				\node[point,red] at (4,3){};
				\node[point,red] at (4,4){};
				\node[point,red] at (4,6){};
				\node[point,red,label={\color{Mycolor2} $\bf {\hspace{-0.25in} n_{\displaystyle \hspace{0.01in}\scalebox{1.075}{$\bf s_{\displaystyle\hspace{0.01in}\scalebox{0.9} {$\bf k$}}$}}}$}] at (4,6){};
				\node[point,red, label={$\bf \color{Mycolor2}\hspace{0.25in}{n_{3}}$}] at (5,0){};
				\node[point,red] at (5,2){};
				\node[point,red] at (5,3){};
				\node[point,red] at (5,5){};
				\node[point,red] at (5,6){};
				\node[point,red] at (6,1){};
				\node[point,red] at (6,2){};
				
				\node[point,red,label={}] at (6,4){};
				\node[point,red] at (6,5){};
				\node[point,red] at (7,0){};
				\node[point,red] at (7,1){};
				\node[point,red, label={}] at (7,3){};
				\node[point,red] at (7,4){};
				\node[point,blue] at (1,0){};
				\node[point,blue] at (1,1){};
				\node[point,blue] at (1,2){};
				\node[point,blue] at (2,1){};
				\node[point,blue] at (2,2){};
				\node[point,blue] at (2,0){};
				\node[point,blue] at (3,0){};
				\node[point,blue] at (3,1){};
				\node[point,blue] at (1,3){};
				\node[point,blue] at (1,4){};
				\node[point,blue] at (2,3){};
				\node[point,blue] at (3,2){};
				\node[point,blue] at (3,4){};
				\node[point,blue] at (3,5){};
				\node[point,blue] at (3,7){};
				\node[point,blue] at (3,8){};
				\node[point,blue] at (2,9){};
				\node[point,blue] at (2,8){};
				\node[point,blue] at (2,5){};
				\node[point,blue] at (2,6){};
				\node[point,blue] at (1,10){};
				\node[point,blue] at (1,9){};
				\node[point,blue] at (1,6){};
				\node[point,blue] at (1,7){};
				\node[point,blue] at (8,0){};
				\node[point,blue] at (8,2){};
				\node[point,blue] at (9,1){};
				\node[point,blue] at (10,0){};
				\node[point,blue] at (11,0){};
				\node[point,blue] at (10,1){};
				\node[point,blue] at (9,2){};
				\node[point,blue] at (8,3){};
				

				
				\draw[dashed,thick,gray](7,4)--node[above=3mm,right]{}(11,0); 
				\draw[dashed,thick,gray](1,9)--node[above=3mm,right]{}(3.3,6.7); 
				\draw[dashed,thick,gray](7,3)--node[above=3mm,right]{}(10,0); 
				\draw[dashed,thick,gray](1,10)--node[above=3mm,right]{}(3.25,7.75); 
				\draw[dashed,thick,gray](1,7)--node[above=3mm,right]{}(3.25,4.75); 
				\draw[dashed,thick,gray](1,6)--node[above=3mm,right]{}(3.25,3.75); 
				\draw[dashed,thick,gray](1,4)--node[above=3mm,right]{}(3.25,1.75); 
				\draw[dashed,thick,gray](1,3)--node[above=3mm,right]{}(3.25,0.75); 
				\draw[dashed,thick,gray](1,2)--node[above=3mm,right]{}(3,0); 
				\draw[dashed,thick,gray](1,1)--node[above=3mm,right]{}(2,0); 
				\draw[dashed,thick,gray](7,1)--node[above=3mm,right]{}(8,0); 

				\node[label={\color{Mycolor2}{\bf {\large I}}}] at (5.5,3.25){};
				
				\node[label={\color{Mycolor2}\scalebox{1.25}{$\bf \displaystyle {{I_{{p}}}}$}}] at (4,-1.25){};
				
				\node[label={\color{Mycolor2}\scalebox{1.25}{$\bf \displaystyle {{I_{{q-1}}}}$}}] at (7,-1.25){};

				\node[label={\color{purple}{\scalebox{1.15}{the line with the equation }}}] at (13.5,1.25){};
				
				\node[label={\color{purple}{\scalebox{1.15}{$x+y=q+k-2$}}}] at (13.35,0.45){};
				\draw[->,thick, forestgreen](9.85,0.25)--node[above=3mm,right]{}(10.75,0.95); 
				\node[label={\color{purple}{\scalebox{1.15}{the line with the equation}}}] at (13.5,3.35){};
				\node[label={\color{purple}{\scalebox{1.15}{$x+y=q+k-1$}}}] at (13.35,2.55){};
				\draw[->,thick, forestgreen](9.65,2)--node[above=3mm,right]{}(11.05,3.3); 
				
		\end{tikzpicture} }
	\end{center}
	\par Now we assume that the terms  $e^{\displaystyle c_{pq}\hspace{0.025in}\alpha_{m}}$, $m\in \mathbb{N}-I$, are on the blue points and the terms $e^{\displaystyle \left(c_{pq}-1\right)\alpha_{ \displaystyle \scalebox{0.9}{$n_{i}$}}}$, $n_{i}\in I$, are on the red points at this line. Before sorting the terms of the sequence $\displaystyle \left(a_{p,n}\over a_{q,n}\right)_{n\in \mathbb{N}}$, we note that the terms of the sequences $\left( e^{\displaystyle c_{pq}\alpha_{m}}\right)_{m\in \mathbb{N}-I}$  and $\left( e^{\displaystyle \left(c_{pq}-1\right)\alpha_{\displaystyle \scalebox{0.9}{$n_{i}$}}}\right) _{i\in \mathbb{N}}$ have decreasing order in themselves.
	\par At first, we take into account the part of $\displaystyle \left(a_{p,n}\over a_{q,n}\right)_{n\in \mathbb{N}}$ including the first $n_{1}-1$ terms $ e^{\displaystyle \left(c_{pq}-1\right)\alpha_{\displaystyle \scalebox{0.9}{$n_{i}$}}}$, $1\leq i\leq n-1$. Since $\alpha$ is increasing, this part has decreasing order and all terms in this part is greater than the terms corresponding to the elements of $I$. Then, having decreasing order, this part remains the same. However, we write this part by shifting to the left taking into account the zero indices for Kolmogorov diameter.
	
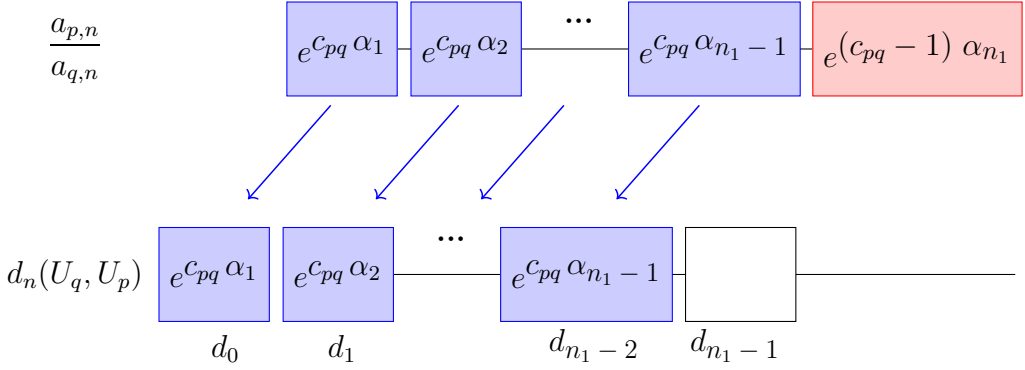
\begin{center}
		\begin{tikzpicture}[point/.style={fill,circle,inner sep=1pt,black}]
			\draw[-](2,-0.5)--node[above=3mm,right]{}(9.2,-0.5);
			\node[label={{$\displaystyle {{a_{p,n}\over a_{q,n}}}$}}] at (-2,-1.25){};
			\node[Myneww Style1] at (1.55,-.5){${e^{\displaystyle  c_{pq}\hspace{0.025in}\alpha_{1}}}$};
			\node[Myneww Style1] at (3.2,-.5){$e^{\displaystyle  c_{pq}\hspace{0.025in}\alpha_{2}}$};
			\node[Myneww Style1] at (6.5,-.5){$e^{\displaystyle  c_{pq}\hspace{0.025in}\alpha_{\displaystyle \scalebox{0.85}{$n_{1}-1$}}}$};
			\node[label={\bf{...}}] at (4.7,-0.45){};
			\node[label={\bf{...}}] at (3,-3.35){};
			\node[Myneww Style2] at (9.2,-.5){$e^{\displaystyle \left( c_{pq}-1\right)\hspace{0.025in}\alpha_{\displaystyle \scalebox{0.85}{$n_{1}$}}}$};
			\draw[-](1,-3.5)--node[above=3mm,right]{}(10.5,-3.5);
			\node[label={{$ d_{n}(U_{q},U_{p})$}}] at (-2,-4.025){};
			\node[Myneww Style1] at (-0.15,-3.5){$e^{\displaystyle  c_{pq}\hspace{0.025in}\alpha_{1}}$};
			\node[Myneww Style1] at (1.5,-3.5){$e^{\displaystyle  c_{pq}\hspace{0.025in}\alpha_{2}}$};
			\node[Myneww Style1] at (4.8,-3.5){$e^{\displaystyle  c_{pq}\hspace{0.025in}\alpha_{\displaystyle \scalebox{0.85}{$n_{1}-1$}}}$};
			\node[Myneww Style3] at (6.85,-3.5){\color{white}{$e^{\displaystyle  c_{pq}\hspace{0.025in}\alpha_{1}}$}};
			\draw[blue,->,line width=0.2mm](1.4,-1.25)--node[above=3mm,right]{}(0.3,-2.5);
			\draw[blue,->,line width=0.2mm](3.1,-1.25)--node[above=3mm,right]{}(2,-2.5);
			\draw[blue,->, line width=0.2mm](6.3,-1.25)--node[above=3mm,right]{}(5.2,-2.5);
			\draw[blue,->, line width=0.2mm](4.5,-1.25)--node[above=3mm,right]{}(3.4,-2.5);
			\node[label={{$ d_{0}$}}] at (0,-4.95){};
			\node[label={{$ d_{1}$}}] at (1.55,-4.95){};
			\node[label={{$ d_{\displaystyle \scalebox{0.85}{$n_{1}-2$}}$}}] at (4.9,-4.95){};
			\node[label={{$ d_{\displaystyle \scalebox{0.85}{$n_{1}-1$}}$}}] at (6.775,-4.95){};
		\end{tikzpicture}
	\end{center}
	So, for every $0\leq n\leq n_{1}-2$,
	$$d_{n}(U_{q}, U_{p})= e^{\displaystyle  c_{pq}\hspace{0.0125in} {\alpha}_{\displaystyle \hspace{0.01in}\scalebox{0.85}{$n+1$}}}.$$
	\par In order to find the diameter $\displaystyle {d_{\displaystyle \scalebox{0.85}{${n_{1}-1}$}}(U_{q}, U_{p})}$, we will compare the term ${ e^{\displaystyle  (c_{pq}-1)\hspace{0.0125in}{\alpha}_{\displaystyle \hspace{0.01in}\scalebox{0.85}{${n_{1}}$}}}}$ with the terms $e^{\displaystyle{c_{pq}\hspace{0.015in}\alpha_{\displaystyle\hspace{-.0in}\scalebox{0.8}{$m$}}}}$, $m\in \mathbb{N}-I$, $m> n_{1}$, and the greatest term gives the diameter $\displaystyle d_{\displaystyle \scalebox{0.85}{${n_{1}-1}$}}(U_{q}, U_{p})$:
	$$e^{\displaystyle {(c_{pq}-1)\alpha_{\hspace{0.01in}\scalebox{0.85}{$n_{1}$}}}}\leq e^{\displaystyle{c_{pq}\hspace{0.015in}\alpha_{\displaystyle\hspace{-.025in}\scalebox{0.8}{ $m$}}}}\hspace{0.25in}\Leftrightarrow\hspace{0.25in} \alpha_{\displaystyle\hspace{-.025in}\scalebox{0.8}{ $m$}}\leq A_{\displaystyle \scalebox{0.8}{$pq$}}\hspace{0.025in}\alpha_{\hspace{0.01in}\scalebox{0.85}{$n_{1}$}}.$$
	where $\displaystyle A_{\displaystyle \scalebox{0.8}{$pq$}}=1+{pq\over {q-p}}$. Then, the terms $e^{\displaystyle c_{pq}\hspace{0.025in}\alpha_{\displaystyle\hspace{-.025in}\scalebox{0.8}{ $m$}}}$, $m\in \mathbb{N}-I$, $m>n_{1}$, satisfying $\alpha_{\displaystyle\hspace{-.025in}\scalebox{0.8}{ $m$}}\leq A_{\displaystyle \scalebox{0.8}{$pq$}}\hspace{0.025in}\alpha_{\hspace{0.01in}\scalebox{0.85}{$n_{1}$}}$ is greater than the term ${ e^{\displaystyle  (c_{pq}-1)\hspace{0.0125in}{\alpha}_{\displaystyle \hspace{0.01in}\scalebox{0.85}{${n_{1}}$}}}}$. So we must write the terms $e^{\displaystyle c_{pq}\hspace{0.025in}\alpha_{\displaystyle\hspace{-.025in}\scalebox{0.8}{ $m$}}}$, $m\in \mathbb{N}-I$, $m>n_{1}$, satisfying $\alpha_{\displaystyle\hspace{-.035in}\scalebox{0.8}{ $m$}}\leq A_{\displaystyle \scalebox{0.8}{$pq$}}\hspace{0.025in}\alpha_{\hspace{0.01in}\scalebox{0.85}{$n_{1}$}}$ before the term ${ e^{\displaystyle  (c_{pq}-1)\hspace{0.0125in}{\alpha}_{\displaystyle \hspace{0.01in}\scalebox{0.85}{${n_{1}}$}}}}$ in decreasing order.
	\par We call the greatest element $m\in \mathbb{N}-I$ satisfying  $\alpha_{\displaystyle\hspace{-.025in}\scalebox{0.8}{ $m$}}\leq A_{\displaystyle \scalebox{0.8}{$pq$}}\hspace{0.025in}\alpha_{\hspace{0.01in}\scalebox{0.85}{$n_{1}$}}$ as $\bf{i_{1}}$. As shown in the following diagram, we can assume that there exists a $k_{1}> 0$ so that the inequality 
	$$\displaystyle n_{\displaystyle \hspace{0.01in} s_{\displaystyle \hspace{0.01in}\scalebox{0.85}{$k_{1}$}} }<i_{1}<n_{\displaystyle \hspace{0.01in} s_{\displaystyle \hspace{0.01in}  \scalebox{0.8}{$(k_{1}+1)$}} }$$
	holds. 
	
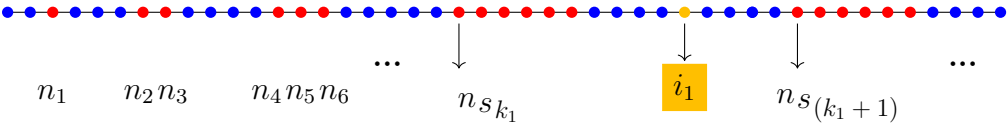
\begin{center}
		\begin{tikzpicture}[point/.style={fill,circle,inner sep=1.5pt,black}]
			
			\draw[-](0,0)--node[above=3mm,right]{}(13.2,0);
			\node[point,blue] at (0,0){};
			\node[point,blue] at (0.3,0){};
			
			\node[point,red] at (0.6,0){};
			\node[label={$n_{1}$}] at (0.6,-1.5){};
			
			\node[point,blue] at (0.9,0){};
			\node[point,blue] at (1.2,0){};
			\node[point,blue] at (1.5,0){};
			
			\node[point,red] at (1.8,0){};
			\node[point,red] at (2.1,0){};
			
			\node[label={$n_{2}$}] at (1.75,-1.5){};
			\node[label={$n_{3}$}] at (2.2,-1.5){};
			
			\node[point,blue] at (2.7,0){};
			\node[point,blue] at (2.4,0){};
			\node[point,blue] at (3,0){};
			\node[point,blue] at (3.3,0){};
			
			\node[point,red] at (3.6,0){};
			\node[point,red] at (3.9,0){};
			\node[point,red] at (4.2,0){};
			
			\node[label={$n_{4}$}] at (3.45,-1.5){};
			\node[label={$n_{5}$}] at (3.9,-1.5){};
			\node[label={$n_{6}$}] at (4.35,-1.5){};

			\node[point,blue] at (4.5,0){};
			\node[point,blue] at (4.8,0){};
			\node[point,blue] at (5.1,0){};
			\node[point,blue] at (5.4,0){};
			\node[point,blue] at (5.7,0){};
			\node[point,red] at (6,0){};
			\node[label={\bf{...}}] at (5.05,-1){};

			\node[point,red] at (6.3,0){};
			\node[point,red] at (6.6,0){};
			\node[point,red] at (6.9,0){};
			\node[point,red] at (7.2,0){};
			\node[point,red] at (7.5,0){};
			\node[point,blue] at (7.8,0){};
			
			\draw[->](6.0,-0.15)--node[above=3mm,right]{}(6.0,-0.75);
			\node[label={{$n_{\displaystyle \hspace{0.01in} s_{\displaystyle \hspace{0.01in}\scalebox{0.8}{$k_{1}$}} }$}}] at (6.4,-1.75){};
			
			\node[point,blue] at (8.1,0){};
			\node[point,blue] at (8.4,0){};
			\node[point,blue] at (8.7,0){};
			
			\node[point,amber] at (9.0,0){};
			\node[My Style3] at (9.0,-1){$i_{1}$};
			\draw[->] (9.0,-0.15)--node[above=3mm,right]{}(9.0,-0.65);
			
			\node[point,blue] at (9.3,0){};
			\node[point,blue] at (9.9,0){};
			\node[point,blue] at (9.6,0){};
			\node[point,blue] at (9.9,0){};
			
			\node[label={{$n_{\displaystyle \hspace{0.01in} s_{\displaystyle \hspace{0.01in}  \scalebox{0.8}{$(k_{1}+1)$}} }$}}] at (11.05,-1.75){};
			\draw[->](10.5,-.15)--node[above=3mm,right]{}(10.5,-0.75);
			
			\node[point,blue] at (10.2,0){};
			\node[point,red] at (10.5,0){};
			\node[point,red] at (10.8,0){};
			\node[point,red] at (11.1,0){};
			\node[point,red] at (11.4,0){};
			\node[point,red] at (11.7,0){};
			\node[point,red] at (12,0){};
			\node[point,blue] at (12.3,0){};
			\node[point,blue] at (12.6,0){};
			\node[point,blue] at (12.9,0){};
			\node[point,blue] at (13.2,0){};
			\node[label={\bf{...}}] at (12.7,-1){};
		\end{tikzpicture}
	\end{center}
	\par This means that the number of elements of $I$ which is less than $i_{1}$ is $\boldsymbol{\displaystyle \scalebox{1.05}{$s_{\displaystyle \hspace{0.01in}  \scalebox{0.75}{$\left(k_{1}+1\right)$}}$}-1}$. So, before the term $ e^{\displaystyle  (c_{pq}-1)\hspace{0.0125in}{\alpha}_{\displaystyle \hspace{0.01in}\scalebox{0.85}{$n_{1}$}}}$, we will write $\displaystyle i_{1}-[\displaystyle\scalebox{1.1}{$s_{\displaystyle \hspace{0.01in}  \scalebox{0.8}{$\left(k_{1}+1\right)$}}$}-1]$ many  $e^{\displaystyle c_{pq}\hspace{0.025in}\alpha_{\displaystyle\hspace{-.025in}\scalebox{0.8}{ $m$}}}$, $m\in \mathbb{N}-I$, $m\leq i_{1}$, terms in decreasing order. Furthermore, while writing these terms in decreasing order, every term ${ e^{\displaystyle  (c_{pq}-1)\hspace{0.0125in}{\alpha}_{\displaystyle \hspace{0.01in}\scalebox{0.85}{${n_{\displaystyle a}}$}}}}$, $1\leq a\leq \displaystyle \scalebox{1.1}{$s_{\displaystyle \hspace{0.01in}  \scalebox{0.8}{$\left(k_{1}+1\right)$}}$}-1$ shifts to the right and every term $e^{\displaystyle c_{pq}\hspace{0.025in}\alpha_{\displaystyle\hspace{-.025in}\scalebox{0.8}{ $m$}}}$, ${m\in \mathbb{N}-I}$, $ m\leq i_{1}$, shifts to the left, as shown in {Diagram  1}.
	\begin{landscape}
		\begin{center}
			\begin{tikzpicture}[point/.style={fill,circle,inner sep=.05pt,black}, rotate=0, transform shape]
				\draw[-](0.65,4)--node[above=3mm,right]{}(21,4);
				\draw[-](0.65,0)--node[above=3mm,right]{}(21,0);
				\draw[-](0.25,-4)--node[above=3mm,right]{}(21,-4);
				\node[label={{$ \hspace{-0.4in}\displaystyle {a_{p,n}\over a_{q,n}}$}}] at (-1.25, 3.325){};
				\node[label={{$\hspace{-0.1in} d_{n}(U_{q}, U_{p})$}}] at (-1.5, -4.525){};
				
				\node[Mynew Style1] at (.65,0){{$1$}};
				\node[Mynew Style1] at (.65, 4){{$1$}};
				\node[Mynew Style1] at (0,-4){{$1$}};
				\draw[->, blue](0.65,-0.5)--node[above=3mm,right]{}(0,-3.5);
				\node[label={{$ \displaystyle d_{0}$}}] at (0,-5.25){};
				
				\node[Mynew Style1] at (1.3,0){{$2$}};
				\node[Mynew Style1] at (1.3,4){{$2$}};
				\node[Mynew Style1] at ( 0.65,-4){{$2$}};
				\draw[->, blue](1.3,-.5)--node[above=3mm,right]{}(0.65,-3.5);
				\node[label={{$ \displaystyle d_{1}$}}] at (0.65,-5.25){};
				
				\node[label={\bf{...}}] at (2.2, 0.015){};
				\node[label={\bf{...}}] at (2.2, 4.015){};
				\node[label={\bf{...}}] at (1.55,-4.15){};
				
				\node[Mynew Style1] at (3.1,0){{}};
				\node[Mynew Style1] at (3.1,4){{}};
				\node[Mynew Style1] at (2.45,-4){{}};
				\draw[->, blue](3.1,-.5)--node[above=3mm,right]{}(2.45,-3.5);

				\node[Mynew Style2] at (3.8,4){{$n_{1}$}};
				\node[Mynew Style3] at (3.8,0){{}};

				\node[Mynew Style1] at (4.5,0){{}};
				\node[Mynew Style1] at (4.5,4){{}};
				\node[Mynew Style1] at (3.1,-4){{}};
				\draw[->, blue](4.5,-0.5)--node[above=3mm,right]{}(3.1,-3.5);
				
				\node[label={\bf{...}}] at (5.4, -0.15){};
				\node[label={\bf{...}}] at (5.4, 4.15){};
				\node[label={\bf{...}}] at (4,-4.15){};
				
				\node[Mynew Style1] at (6.3,0){{}};
				\node[Mynew Style1] at (6.3,4){{}};
				\node[Mynew Style1] at (4.9,-4){{}};
				\draw[->, blue](6.3,-0.5)--node[above=3mm,right]{}(4.9,-3.5);

				\node[Mynew Style2] at (7,4){{$n_{2}$}};
				\node[Mynew Style2] at (7.75,4){{$n_{3}$}};
				\node[Mynew Style3] at (7,0){{}};
				\node[Mynew Style3] at (7.75,0){{}};
				
				\node[Mynew Style1] at (8.45,0){{}};
				\node[Mynew Style1] at (8.45,4){{}};
				\node[Mynew Style1] at (5.6,-4){{}};
				\draw[->, blue](8.45,-0.5)--node[above=3mm,right]{}(5.6,-3.5);
				
				\node[label={\bf{...}}] at (9.35, -0.15){};
				\node[label={\bf{...}}] at (9.35, 4.15){};
				\node[label={\bf{...}}] at (6.5,-4.15){};
				
				\node[Mynew Style1] at (10.25,0){{}};
				\node[Mynew Style1] at (10.25,4){{}};
				\node[Mynew Style1] at (7.4,-4){{}};
				\draw[->, blue](10.25,-0.5)--node[above=3mm,right]{}(7.4,-3.5);

				\node[Mynew Style2] at (10.95,4){{$n_{4}$}};
				\node[Mynew Style2] at (11.675,4){{$n_{5}$}};
				\node[Mynew Style2] at (12.4,4){{$n_{6}$}};
				\node[Mynew Style3] at (10.95,0){{}};
				\node[Mynew Style3] at (11.675,0){{}};
				\node[Mynew Style3] at (12.4,0){{}};
				
				\node[Mynew Style1] at (13.1,0){{}};
				\node[Mynew Style1] at (13.1,4){{}};
				\node[Mynew Style1] at (8.1,-4){{}};
				\draw[->, blue](13.1,-0.5)--node[above=3mm,right]{}(8.1,-3.5);
				
				\node[label={\bf{...}}] at (14, -0.15){};
				\node[label={\bf{...}}] at (14, 4.15){};
				\node[label={\bf{...}}] at (9,-4.15){};
				
				\node[Mynew Style1] at (14.9,0){{}};
				\node[Mynew Style1] at (14.9,4){{}};
				\node[Mynew Style1] at (9.9,-4){{}};
				\draw[->, blue](14.9,-0.5)--node[above=3mm,right]{}(9.9,-3.5);

				\node[label={{${n_{\displaystyle \hspace{0.01in} \scalebox{1.1}{$\displaystyle \hspace{0.01in} s_{\displaystyle \hspace{0.01in} \scalebox{0.8}{$k_{\displaystyle \scalebox{0.8}{1}}$}}$}}}$}}] at (15.55, 5){};
				\draw[->](15.5, 4.95)--node[above=3mm,right]{}(15.5, 4.5);
				\node[Mynew Style2] at (15.55,4){{}};11
				\node[label={\bf{...}}] at (16.45, 4.15){};
				\node[Mynew Style2] at (17.35,4){{}};
				
				\node[Mynew Style3] at (15.55,0){{}};
				\node[label={\bf{...}}] at (16.45, -0.15){};
				\node[Mynew Style3] at (17.35,0){{}};
				
				\node[Mynew Style1] at (18.05,0){{}};
				\node[Mynew Style1] at (18.05,4){{}};
				\node[Mynew Style1] at (10.6,-4){{}};
				\draw[->, blue](18.05,-0.5)--node[above=3mm,right]{}(10.6,-3.5);
				
				\node[label={\bf{...}}] at (18.95, -0.15){};
				\node[label={\bf{...}}] at (18.95, 4.15){};
				\node[label={\bf{...}}] at (11.5,-4.15){};
				
				\node[Mynew Style1] at (19.85,0){{}};
				\node[Mynew Style1] at (19.85,4){{}};
				\node[Mynew Style1] at (12.4,-4){{}};
				\draw[->, blue](19.85,-0.5)--node[above=3mm,right]{}(12.4,-3.5);
				
				\node[Mynew Style4] at (20.55,0){{$i_{1}$}};
				\node[Mynew Style4] at (20.55,4){{$i_{1}$}};
				\node[Mynew Style4] at (13.1,-4){{$i_{1}$}};
				\draw[->, blue](20.55,-0.5)--node[above=3mm,right]{}(13.1,-3.5);
				
				\node[Mynew Style2] at (13.8,-4){{$n_{1}$}};
				\draw[->, red](3.8, 3.75)--node[above=3mm,right]{}(13.8,-3.5);
				\node[label={{$ \displaystyle d_{\displaystyle \scalebox{0.85}{$j_{\displaystyle \scalebox{0.8}{$1$}}$}}$}}] at (13.8,-5.25){};
				\draw[->,dashed, red](7, 3.6)--node[above=3mm,right]{}(17.6,-3.5);
				\draw[->,dashed, red](7.75, 3.6)--node[above=3mm,right]{}(18,-3.5);
				\draw[->,dashed, red](10.95, 3.6)--node[above=3mm,right]{}(18.4,-3.5);
				\draw[->,dashed, red](11.675, 3.6)--node[above=3mm,right]{}(18.8,-3.5);
				\draw[->,dashed, red](12.4, 3.6)--node[above=3mm,right]{}(19.2,-3.5);
				\draw[->,dashed, red](15.55, 3.6)--node[above=3mm,right]{}(19.6,-3.5);
				\draw[->,dashed, red](17.35, 3.6)--node[above=3mm,right]{}(20,-3.5);
				
				\node[label={\textnormal{\bf Diagram 1}}] at (10, -7.0){};
				\node[label={\textnormal{ For decreasing order, we will lift the red boxes and shift the blue boxes to the left.}}] at (10, 1.5){};
			\end{tikzpicture}
		\end{center}
	\end{landscape}
	\par In order to find  $n_{1}-1$-th Kolmogorov diameter, we shift the term corresponding to the first element $n_{1}$ of $I$. Considering also that we shift the terms to the left for $d_{0}(U_{q}, U_{p})$, we find that for every  $n_{1}-1\leq n\leq n_{2}-3$,
	$${ \displaystyle d_{n}(U_{q}, U_{p})=e^{\mathlarger{c_{pq}\hspace{0.015in}{\alpha}_{\displaystyle \hspace{0.01in} \scalebox{0.85}{$n+2$}}}}}.$$
	So, we found the Kolmogorov diameters until the indices $n_{2}-2$. Now, we also shift the terms corresponding to the element $n_{2}$ and $n_{3}$ of $I$. Up till now, we shift the terms to the left four-indices, then we find that for every ${n_{2}-2\leq n\leq n_{4}-5}$
	$${\displaystyle d_{n}(U_{q}, U_{p})=e^{\mathlarger{c_{pq}\hspace{0.015in}{\alpha}_{\hspace{0.01in} \scalebox{0.85}{$n+4$}}}}}.$$
	\par We would like to point out that the endpoints of the intervals in which we determine Kolmogorov diameters are generally represented by the elements of $I_{p}$. Because the terms corresponding to the elements of $I$ that we shift to the right and the terms corresponding to the elements of $\mathbb{N}-I$ that we shift to the left are between the two elements of $I_{p}$, as seen in the following diagram.
	\begin{center}\scalebox{0.9}
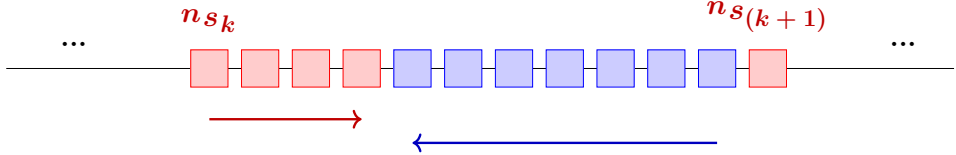
{
			\begin{tikzpicture}[point/.style={fill,circle,inner sep=.05pt,black}, rotate=0, transform shape]
				\draw[-](-2.0,4)--node[above=3mm,right]{}(12,4);

				\node[Mynew Style2] at (1, 4){};
				\node[Mynew Style2] at (1.75, 4){};
				\node[Mynew Style2] at (2.5, 4){};
				\node[Mynew Style2] at (3.25, 4){};
				
				\draw[->, line width=0.35mm, Mycolor2](1, 3.25)--node[above=3mm,right]{}(3.25, 3.25);

				\draw[<-, line width=0.35mm, Mycolor4](4, 2.9)--node[above=3mm,right]{}(8.5, 2.9);

				\node[Mynew Style1] at (4,4){};
				\node[Mynew Style1] at (4.75,4){};
				\node[Mynew Style1] at (5.5,4){};
				\node[Mynew Style1] at (6.25,4){};
				\node[Mynew Style1] at (7, 4){};
				\node[Mynew Style1] at (7.75,4){};
				\node[Mynew Style1] at (8.5, 4){};
				
				\node[Mynew Style2] at (9.25, 4){};
				
				\node[label={\bf{...}}] at (-1, 4.05){};
				\node[label={\bf{...}}] at (11.25, 4.05){};
				
				\node[label={\color{Mycolor2}{$\boldsymbol{n_{\displaystyle \hspace{0.01in} \scalebox{1.05}{$s_{\displaystyle \scalebox{0.8}{$\boldsymbol{(k+1)}$}} $}}}$}}] at (9.25,4.25){};
				
				\node[label={\color{Mycolor2}{$\boldsymbol{ n_{\displaystyle \hspace{0.01in} \scalebox{1.05}{$s_{\displaystyle \hspace{0.01in}\scalebox{0.8}{$\boldsymbol{k}$}} $}}}$}}]at (1,4.25){};
				
			\end{tikzpicture}
		}
	\end{center}
	Another significant point in writing the endpoints of the intervals in which we determine the diameters are to find out how many elements, the terms corresponding to the elements of $I$, we shift to the right.  
	\par We continue to calculate the diameters with this perspective. Let us assume that we replaced $n_{\displaystyle \hspace{0.01in} \scalebox{1.05}{$s_{\displaystyle \scalebox{0.65}{$0$}}$}}-[\scalebox{1.05}{$s_{\displaystyle \scalebox{0.65}{$0$}}$}+1]$ terms in decreasing order. In order to find  $\displaystyle n_{\displaystyle \hspace{0.01in} \scalebox{1.05}{$s_{\displaystyle \scalebox{0.65}{$0$}}$}}-\scalebox{1.05}{$s_{\displaystyle \scalebox{0.65}{$0$}}$}$-th Kolmogorov diameter, we shift $\scalebox{1.1}{$s_{1}$}$ terms corresponding to the elements of I in total, for every ${\displaystyle n_{\displaystyle \hspace{0.01in} {\scalebox{1.05}{$s_{\displaystyle \scalebox{0.65}{$0$}}$}}}-\scalebox{1.05}{$s_{\displaystyle \scalebox{0.65}{$0$}}$}\leq n\leq n_{\displaystyle \hspace{0.01in} \scalebox{1.05}{$s_{\displaystyle \scalebox{0.65}{$1$}}$}}-[\scalebox{1.05}{$s_{\displaystyle \scalebox{0.65}{$1$}}$}+1]}$, we have
	$${\displaystyle d_{n}(U_{q}, U_{p})=e^{\displaystyle {c_{pq}\hspace{0.015in} {\alpha}_{\displaystyle \hspace{0.01in}\scalebox{0.9}{${n+\scalebox{1.1}{$s_{\displaystyle \scalebox{0.65}{$1$}}$}}$} }}}}.$$
	Considering the terms that we shift to the right in each step, we can write for every $0\leq k< k_{1}$ and for every 
	${n_{\displaystyle \hspace{0.01in} \scalebox{1.05}{$s_{\displaystyle\scalebox{0.75} {$k$}}$}} - \scalebox{1.05}{$s_{\displaystyle\scalebox{0.8} {$k$}}$}\leq n\leq n_{\displaystyle \hspace{0.01in} \scalebox{1.05}{$s_{\displaystyle\scalebox{0.75} {$(k+1)$}}$}} -[\scalebox{1.05}{$s_{\displaystyle\scalebox{0.75} {$(k+1)$}}$}+1]}$
	$${\displaystyle d_{n}\left(U_{q},U_{p}\right)=e^{\displaystyle c_{pq}\hspace{0.025in}\alpha_{\displaystyle\hspace{-.025in}\scalebox{0.9}{ $n+\scalebox{1.1}{$s_{\displaystyle\scalebox{0.8} {$(k+1)$}}$}$}}}} $$
	and for all ${n_{\displaystyle \hspace{0.01in} \scalebox{1.05}{$s_{\displaystyle \hspace{0.01in}\scalebox{0.8}{$k_{1}$}}$}}-\scalebox{1.05}{$s_{\displaystyle \hspace{0.01in}\scalebox{0.8}{$k_{1}$}}$} \leq n\leq i_{1}-\scalebox{1.05}{$s_{\displaystyle \hspace{0.01in}\scalebox{0.75}{$(k_{1}+1)$}}$}}$
	$${\displaystyle d_{n}\left(U_{q},U_{p}\right)=e^{\displaystyle c_{pq}\hspace{0.025in}\alpha_{\displaystyle\hspace{-.025in}\scalebox{0.9}{ $n+\scalebox{1.1}{$s_{\displaystyle\scalebox{0.8} {$(k_{1}+1)$}}$}$}}}}.$$
	Therefore, we shift $\displaystyle i_{1}-[ \scalebox{1.05}{$s_{\displaystyle \hspace{0.01in}\scalebox{0.8}{$(k_{1}+1)$}}$}-1]$ many terms $e^{\displaystyle c_{pq}\hspace{0.025in}\alpha_{\displaystyle\hspace{-.025in}\scalebox{0.8}{ $m$}}}$, $m\in \mathbb{N}-I$, $m\leq i_{1}$ to left, namely, we sort all terms which is greater than ${ e^{\displaystyle  (c_{pq}-1)\hspace{0.0125in}{\alpha}_{\displaystyle \hspace{0.01in}\scalebox{0.85}{${n_{1}}$}}}}$. Hence, the term ${ e^{\displaystyle  (c_{pq}-1)\hspace{0.0125in}{\alpha}_{\displaystyle \hspace{0.01in}\scalebox{0.85}{${n_{1}}$}}}}$
	is replaced at the indices $ {\boldsymbol{ j_{1}}=i_{1}-\scalebox{1.05}{$s_{\displaystyle \hspace{0.01in}\scalebox{0.8}{$(k_{1}+1)$}}$}+1}$, namely,
	$${ d_{\displaystyle \hspace{0.01in}\scalebox{0.85}{$j_{1}$}}}(U_{q}, U_{p})={ e^{\displaystyle  (c_{pq}-1)\hspace{0.0125in}{\alpha}_{\displaystyle \hspace{0.01in}\scalebox{0.88}{${n_{1}}$}}}}.$$
\par Now assume that the first $a-1$, $(a\geq 2)$ terms corresponding to the elements of $I$ are placed in decreasing order. Before the term ${ e^{\displaystyle  (c_{pq}-1)\hspace{0.0125in}{\alpha}_{\displaystyle \hspace{0.01in}\scalebox{0.85}{${n_{\displaystyle \scalebox{0.85}{$a$}}}$}}}}$, we must write the terms $e^{\displaystyle c_{pq}\hspace{0.025in}\alpha_{\displaystyle\hspace{-.025in}\scalebox{0.8}{ $m$}}}$, $m\in \mathbb{N}-I$ which is greater than ${ e^{\displaystyle  (c_{pq}-1)\hspace{0.0125in}{\alpha}_{\displaystyle \hspace{0.01in}\scalebox{0.85}{${n_{\displaystyle \scalebox{0.85}{$a$}}}$}}}}$, satisfying
the inequality
$\alpha_{\displaystyle\hspace{-.025in}\scalebox{0.8}{ $m$}}\leq A_{\displaystyle \scalebox{0.8}{$pq$}}\hspace{0.025in}\alpha_{\hspace{0.01in}\scalebox{0.85}{${n_{\displaystyle \scalebox{0.85}{$a$}}}$}}$ We call the greatest element of $m\in \mathbb{N}$ satisfying 
$\alpha_{\displaystyle\hspace{-.025in}\scalebox{0.8}{ $m$}}\leq A_{\displaystyle \scalebox{0.8}{$pq$}}\hspace{0.025in}\alpha_{\hspace{0.01in}\scalebox{0.85}{${n_{\displaystyle \scalebox{0.85}{$a$}}}$}}$ as $\boldsymbol{i_{\scalebox{0.85}{$a$}}}$. We can assume that there exists a $k_{\scalebox{0.85}{$a$}}\in \mathbb{N}$ so that
$$\displaystyle n_{\displaystyle \hspace{0.01in} \scalebox{1.05}{$s_{\displaystyle \hspace{0.01in}\scalebox{0.8}{$k_{\hspace{0.01in}\scalebox{0.9}{$a$}}$}} $}}<i_{\hspace{0.01in} \scalebox{0.85}{$a$}}<n_{\displaystyle \hspace{0.01in} \scalebox{1.05}{$s_{\displaystyle \hspace{0.01in}  \scalebox{0.8}{$(k_{\hspace{0.01in}\scalebox{0.9}{$a$}}+1)$}}$} }.$$ 
This means that the number of elements of $I$ which is less than $i_{\scalebox{0.85}{$a$}}$ is ${\scalebox{1.05}{$s_{\displaystyle \hspace{0.01in}  \scalebox{0.8}{$(k_{\hspace{0.01in}\scalebox{0.9}{$a$}}+1)$}}$}-1}$. So, before the term ${e^{\displaystyle  (c_{pq}-1)\hspace{0.0125in}{\alpha}_{\displaystyle \hspace{0.01in}\scalebox{0.85}{${n_{\displaystyle \scalebox{0.85}{$a$}}}$}}}}$, we write $\displaystyle i_{\scalebox{0.85}{$a$}}-(\scalebox{1.05}{$s_{\displaystyle \hspace{0.01in} \scalebox{0.8}{$(k_{ \scalebox{0.85}{$a$}}+1)$}}$}-1)$ many $e^{\displaystyle c_{pq}\hspace{0.025in}\alpha_{\displaystyle\hspace{-.025in}\scalebox{0.8}{ $m$}}}$, $m\in \mathbb{N}-I$, $m\leq i_{\scalebox{0.85}{$a$}}$, terms in decreasing order. Since we assume that the first $a-1$ terms corresponding to the elements of $I$ is placed in decreasing order, then the term ${ e^{\displaystyle  (c_{pq}-1)\hspace{0.0125in}{\alpha}_{\displaystyle \hspace{0.01in}\scalebox{0.85}{${n_{\displaystyle \scalebox{0.85}{$a$}}}$}}}}$ are replaced at the indices ${\boldsymbol{ j_{\hspace{0.01in}\scalebox{0.85}{$a$}}}}=i_{\hspace{0.01in}\scalebox{0.85}{$a$}}-\scalebox{1.05}{$s_{\displaystyle \hspace{0.01in}  \scalebox{0.85}{$(k_{\scalebox{0.9}{$\hspace{0.01in}a$}}+1)$}}$}-1+a$, namely,
$${d_{\displaystyle \scalebox{0.85}{$j_{\scalebox{0.85}{\hspace{0.01in}$a$}}$}}(U_{q}, U_{p})={ e^{\displaystyle  (c_{pq}-1)\hspace{0.0125in}{\alpha}_{\displaystyle \hspace{0.01in}\scalebox{0.85}{${n_{\displaystyle \scalebox{0.9}{\hspace{0.01in}$a$}}}$}}}}.}$$
\par Now, we determine Kolmogorov diameters between the indices $j_{\scalebox{0.85}{\hspace{0.01in}$a$}}$ and $j_{\scalebox{0.85}{\hspace{0.01in}$(a+1)$}}$ for every $a\geq 1$. Starting the index $j_{\scalebox{0.85}{$a$}}$, we must compare the term ${ e^{\displaystyle  (c_{pq}-1)\hspace{0.0125in}{\alpha}_{\displaystyle \hspace{0.01in}\scalebox{0.85}{${n_{\hspace{0.01in}\displaystyle \scalebox{0.9}{$a+1$}}}$}}}}$ with  the terms $e^{\displaystyle c_{pq}\hspace{0.025in}\alpha_{\displaystyle\hspace{-.025in}\scalebox{0.75}{ $m$}}}$, $ m\in \mathbb{N}-I$ and for every $ m\in \mathbb{N}-I$ satisfying $\alpha_{\displaystyle\hspace{-.025in}\scalebox{0.8}{ $m$}}\leq A_{\displaystyle \scalebox{0.8}{$pq$}}\hspace{0.025in}\alpha_{\hspace{0.01in}\scalebox{0.85}{${n_{\displaystyle \scalebox{0.9}{$a+1$}}}$}}$, we write the terms $e^{\displaystyle c_{pq}\hspace{0.025in}\alpha_{\displaystyle\hspace{-.025in}\scalebox{0.8}{ $m$}}}$, before the term  ${ e^{\displaystyle  (c_{pq}-1)\hspace{0.0125in}{\alpha}_{\displaystyle \hspace{0.01in}\scalebox{0.85}{${n_{\displaystyle \scalebox{0.9}{$a+1$}}}$}}}}$. Again, we call the largest element of $\mathbb{N}-I$ satisfying above inequality as $i_{\displaystyle \scalebox{0.85}{$(a+1)$}}$ for which there is $k_{\displaystyle \scalebox{0.9}{$(a+1)$}}\in \mathbb{N}$ satisfying
$$\displaystyle n_{\displaystyle \hspace{0.01in} \scalebox{1.05}{$ s_{\displaystyle \hspace{0.01in}\scalebox{0.85}{$k_{\displaystyle\scalebox{0.9}{$(a+1)$}}$}}$} }<i_{\displaystyle \scalebox{0.85}{$(a+1)$}}<n_{\displaystyle \hspace{0.01in} \scalebox{1.05}{$s_{\displaystyle \hspace{0.01in}  \scalebox{0.85}{$(k_{\displaystyle \scalebox{0.9}{$(a+1)$}}+1)$}} $}}.$$
\par Let us continue to decreasing order from $j_{\scalebox{0.9}{$a$}}+1$:~\\For all ${j_{\scalebox{0.85}{$a$}}+1\leq n\leq 
	n_{\displaystyle \hspace{0.01in} \scalebox{1.05}{$s_{\displaystyle \hspace{0.01in}\scalebox{0.9}{$(k_{\displaystyle\scalebox{0.9}{$a$}}+1)$}}$} }- \scalebox{1.05}{$s_{\displaystyle \hspace{0.01in}\scalebox{0.85}{$(k_{\displaystyle\scalebox{0.9}{$a$}}+1)$}}$}+a-1}$ 
$${ d_{n}\left(U_{q},U_{p}\right)= e^{\displaystyle {c_{pq}\hspace{0.015in}\alpha_{\displaystyle\hspace{0.01in} \scalebox{0.9}{$n+\scalebox{1.1}{$s_{\displaystyle \scalebox{0.85}{$(k_{\scalebox{0.9}{$a$}}+1)$}}$}-a$}}}}.}$$
If any, for every $k_{\displaystyle \scalebox{0.85}{$a$}}+1\leq k\leq k_{\displaystyle \scalebox{0.85}{$(a+1)$}}-1$ and for every  ~\\ ${n_{\displaystyle \hspace{0.01in} \scalebox{1.05}{$s_{\displaystyle \hspace{0.01in}\scalebox{0.8}{$k$}}$}}-\scalebox{1.05}{$s_{\displaystyle \hspace{0.01in}\scalebox{0.85}{$k$}}$}+a\leq n\leq n_{\displaystyle \hspace{0.01in} \scalebox{1.05}{$s_{\displaystyle \hspace{0.01in}\scalebox{0.8}{$(k+1)$}}$}}-\scalebox{1.05}{$s_{\displaystyle \hspace{0.01in}\scalebox{0.8}{$(k+1)$}}$}+a-1}$
$${d_{n}\left(U_{q},U_{p}\right)= e^{\displaystyle {c_{pq}\hspace{0.015in}\alpha_{\displaystyle \scalebox{0.9}{$n+\scalebox{1.05}{$s_{\displaystyle \scalebox{0.9}{$(k+1)$}}$}-a$}}}}}$$
and for every ${ n_{\displaystyle \hspace{0.01in} \scalebox{1.05}{$s_{\displaystyle \hspace{0.01in}  \scalebox{0.85}{$k_{\displaystyle \scalebox{0.9}{$(a+1)$}}$}}$} } -\scalebox{1.05}{$s_{\displaystyle \hspace{0.01in}  \scalebox{0.85}{$k_{\scalebox{0.9}{$(a+1)$}}$}}$}+a\leq n\leq i_{\scalebox{0.85}{$(a+1)$}}-\scalebox{1.05}{$s_{\displaystyle \hspace{0.01in}  \scalebox{0.85}{$(k_{\scalebox{0.9}{$(a+1)$}}+1)$}}$}+a}$
$${d_{n}\left(U_{q},U_{p}\right)=e^{\displaystyle {c_{pq}\hspace{0.015in}\alpha_{\displaystyle \scalebox{0.9}{$n+\scalebox{1.1}{$s_{\displaystyle \hspace{0.01in}  \scalebox{0.85}{$(k_{\scalebox{0.9}{$(a+1)$}}+1)$}}$}-a$}}}}}.$$
\par We sort all terms which is greater than ${ e^{\displaystyle  (c_{pq}-1)\hspace{0.0125in}{\alpha}_{\displaystyle \hspace{0.01in}\scalebox{0.85}{${n_{\displaystyle \scalebox{0.9}{$a+1$}}}$}}}}$. Then, the term ${ e^{\displaystyle  (c_{pq}-1)\hspace{0.0125in}{\alpha}_{\displaystyle \hspace{0.01in}\scalebox{0.85}{${n_{\displaystyle \scalebox{0.9}{$a+1$}}}$}}}}$ is replaced at the indices
\\ ${ j_{\scalebox{0.85}{$(a+1)$}}=i_{\scalebox{0.85}{$(a+1)$}}-\scalebox{1.05}{$s_{\displaystyle \hspace{0.01in}  \scalebox{0.85}{$(k_{\scalebox{0.9}{$(a+1)$}}+1)$}}$}+a+1}$, namely,
$${d_{\displaystyle \scalebox{0.9}{$j_{\scalebox{0.85}{$(a+1)$}}$}}(U_{q}, U_{p})={ e^{\displaystyle  (c_{pq}-1)\hspace{0.0125in}{\alpha}_{\displaystyle \hspace{0.01in}\scalebox{0.85}{${n_{\displaystyle \scalebox{0.9}{$(a+1)$}}}$}}}}.}$$
Hence, we determine all Kolmogorov diameters between the terms 
${e^{\displaystyle  (c_{pq}-1)\hspace{0.0125in}{\alpha}_{\displaystyle \hspace{0.01in}\scalebox{0.9}{${n_{\displaystyle \scalebox{0.85}{$a$}}}$}}}}$
and  ${ e^{\displaystyle  (c_{pq}-1)\hspace{0.0125in}{\alpha}_{\displaystyle \hspace{0.01in}\scalebox{0.85}{${n_{\displaystyle \scalebox{0.9}{$(a+1)$}}}$}}}}$ for every $a\geq 1$. 
\par Therefore, we can calculate all Kolmogorov diameters by following the above observation, and finally we can write:
~\\ 1. Let $\displaystyle J:=\left\{j_{\scalebox{0.85}{$a$}}: a\in \mathbb{N}\right\}$ where $j_{\scalebox{0.85}{$a$}}=i_{\scalebox{0.85}{$a$}}-\scalebox{1.05}{$s_{\displaystyle \hspace{0.01in}  \scalebox{0.85}{$(k_{a}+1)$}}$}-1+a$. For all $a\in \mathbb{N}$,
$$\boxed{d_{\displaystyle \scalebox{0.85}{$j_{\scalebox{0.85}{$a$}}$}}(U_{q}, U_{p})={ e^{\displaystyle  (c_{pq}-1)\hspace{0.0125in}{\alpha}_{\displaystyle \hspace{0.01in}\scalebox{0.85}{${n_{\displaystyle \scalebox{0.85}{$a$}}}$}}}}.}$$
2. For $a, k\in \mathbb{N}$, we define 
\begin{equation*}
	I_{\displaystyle \scalebox{0.8}{$a,k$}}=\left[n_{\displaystyle \hspace{0.01in} \scalebox{1.005}{$s_{\displaystyle \hspace{0.01in}\scalebox{0.8}{$k$}}$} }-\scalebox{1.005}{$s_{\displaystyle \hspace{0.01in}\scalebox{0.8}{$k$}}$}+a , n_{\displaystyle \hspace{0.01in} \scalebox{1.005}{$s_{\displaystyle \hspace{0.01in}\scalebox{0.8}{$(k+1)$}}$}}-\scalebox{1.005}{$s_{\displaystyle \hspace{0.01in}\scalebox{0.8}{$(k+1)$}}$}+a-1\right]
\end{equation*}
and 
$$\displaystyle K=\bigcup_{\scalebox{0.8}{$a\in \mathbb{N}$}}\hspace{0.1in}\bigcup_{\scalebox{0.8}{$k_{\displaystyle \scalebox{0.9}{$a$}}+1\leq k\leq k_{\displaystyle \scalebox{0.9}{$a+1$}}-1$}} I_{\displaystyle \scalebox{0.8}{$a,k$}}.$$
For every $n\in K$, there is an $a\in \mathbb{N}$ and a $k\in \mathbb{N}$ satisfying $k_{\displaystyle \scalebox{0.8}{$a$}}+1\leq k\leq k_{\displaystyle \scalebox{0.8}{$(a+1)$}}-1$ such that 
$$\boxed{d_{n}\left(U_{q},U_{p}\right)= e^{\displaystyle {c_{pq}\hspace{0.015in}\alpha_{\displaystyle\hspace{0.01in} \scalebox{0.9}{$n+\scalebox{1.05}{$s_{\displaystyle \scalebox{0.8}{$(k+1)$}}$}-a$}}}}.}$$
3. Let $
\displaystyle L= \bigcup_{a\in \mathbb{N}} 
\left[ \right.$ $
j_{\scalebox{0.85}{$a$}}+1,$ $
n_{ \displaystyle \hspace{-0.025in} \scalebox{1.05}{
	${s_{\displaystyle \scalebox{0.85}{$( k_{\scalebox{0.9}{$a$}}+1)$}}}$}}$ $-\scalebox{1.05}{$ s_{\displaystyle \hspace{-0.025in} \scalebox{0.85}{ $(
		k_{\scalebox{0.9}{$a$}}+1 )$}}$}$+$a-1 
\left. \right]$. For every $n\in L$, there is an $a\in \mathbb{N}$ such that
$$\boxed{d_{n}\left(U_{q},U_{p}\right)= e^{\displaystyle {c_{pq}\hspace{0.015in}\alpha_{\displaystyle \scalebox{0.9}{$\hspace{0.01in}n+\scalebox{1.05}{$s_{\displaystyle \scalebox{0.85}{$(k_{\scalebox{0.9}{$a$}}+1)$}}$}-a$}}}}.}$$
4.  Let $\displaystyle M= \bigcup_{a\in \mathbb{N}} \left[ \displaystyle n_{\displaystyle \hspace{0.01in} \scalebox{1.05}{$s_{\displaystyle \hspace{0.0in}  \scalebox{0.85}{$k_{\displaystyle \scalebox{0.85}{$a$}}$}}$} } -\scalebox{1.05}{$s_{\displaystyle \hspace{0.0in}  \scalebox{0.85}{$k_{\scalebox{0.85}{$a$}}$}}$}+a-1, j_{\scalebox{0.85}{$a$}}-1\right]$. For every $n\in M$, there is an $a\in \mathbb{N}$ such that
$$\boxed{d_{n}\left(U_{q},U_{p}\right)=e^{\displaystyle {c_{pq}\hspace{0.015in}\alpha_{\displaystyle \scalebox{0.9}{$n+\scalebox{1.05}{$s_{\displaystyle \hspace{0.01in}  \scalebox{0.85}{$(k_{\scalebox{0.9}{$a$}}+1)$}}$}-(a-1)$}}}}}.$$	
\par All Kolmogorov diameters in the light of above observation are found since $\displaystyle \mathbb{N} = \left\lbrace0, 1,..., n_{1}-2\right\rbrace  \cup J\cup K\cup L\cup M$. This completes the determination of the diameters.
\par Now, we give an estimation for Kolmogorov diameters of an element $\bf \mathcal{K}_{\alpha}$ of the family $\bf \mathcal{K}$ which is parameterized by $\alpha$. 
\begin{thm}\label{stabilitycase} Let $\bf \mathcal{K}_{\alpha}$ be an element of the family $\bf \mathcal{K}$ with the parameter $\alpha$.  For every $p$, $q>p$ there exists a $N\in \mathbb{N}$ such that
\begin{equation} \label{stability} \displaystyle
	e^{\mathlarger{c_{pq}\hspace{0.015in}\alpha_{\displaystyle\hspace{0.01in} \scalebox{0.8}{$4n$}}}}\leq d_{n}(U_{q}, U_{p})\leq e^{\mathlarger{c_{pq}\hspace{0.015in} \alpha_{\displaystyle \scalebox{0.8}{$n$}}}}
\end{equation}
for every $n\geq N$.
\end{thm}
\begin{proof}
Let $p\in \mathbb{N}$ and $q>p$. Above we obtained Kolmogorov diameters $d_{n}(U_{q}, U_{p})$ on each subsets $\left\lbrace 0, 1,..., n_{1}-1\right\rbrace$, $J$, $K$, $L$ and $M$ of $\mathbb{N}$. We will show that the inequality \ref{stability} holds for sufficiently large elements of each subsets  $J$, $K$, $L$, and $M$ of $\mathbb{N}$.
\par Primarily, we will show that ${2j_{\scalebox{0.8}{$\hspace{0.01in}a$}}>i_{\scalebox{0.8}{$\hspace{0.01in}a$}}}$ for sufficiently large $a\in \mathbb{N}$. We know that for every $i_{\scalebox{0.8}{$a$}}$ there exists a $k_{\scalebox{0.8}{$a$}}\in \mathbb{N}$ satisfying
$\displaystyle n_{\displaystyle \hspace{0.01in} \scalebox{1.05}{$s_{\displaystyle \hspace{0.01in}\scalebox{0.85}{$k_{\scalebox{0.9}{$\hspace{0.01in}a$}}$}}$}}<i_{\scalebox{0.9}{$\hspace{0.01in}a$}}<n_{\displaystyle \hspace{0.01in} \scalebox{1.05}{$s_{\displaystyle \hspace{0.01in}  \scalebox{0.85}{$(k_{\scalebox{0.9}{$\hspace{0.01in}a$}}+1)$}}$} }.$ Since $n_{\displaystyle \hspace{0.01in} \scalebox{1.05}{$s_{\displaystyle \hspace{0.01in}\scalebox{0.85}{$k_{\displaystyle \scalebox{0.9}{$\hspace{0.01in}a$}}$}}$}}$ is on the line which has the equation $x+y=q+k_{\scalebox{0.85}{$\hspace{0.01in}a$}}-2$, the first element  of $I_{\displaystyle \scalebox{0.8}{$q+k_{\scalebox{0.9}{$\hspace{0.01in}a$}}-2$}}$ is less than $n_{\displaystyle \hspace{0.01in} \scalebox{1.05}{$s_{\displaystyle \hspace{0.01in}\scalebox{0.85}{$k_{\scalebox{0.9}{$\hspace{0.01in} a$}}$}}$} }$, we can write
$$n_{\displaystyle \hspace{0.01in} \scalebox{1.05}{$s_{\displaystyle \hspace{0.01in}\scalebox{0.85}{$k_{\displaystyle \scalebox{0.85}{$\hspace{0.01in}a$}}$}}$}} \geq {\left(q+k_{\displaystyle \scalebox{0.85}{$\hspace{0.01in}a$}}-2\right)\left(q+k_{\displaystyle \scalebox{0.85}{$\hspace{0.01in}a$}}-1\right)\over 2}= {{\left( q-2\right) \left(q-1\right)}\over 2}+\left( q-1\right) k_{\displaystyle \scalebox{0.85}{$\hspace{0.01in}a$}}+{k_{\displaystyle \scalebox{0.85}{$\hspace{0.01in}a$}}\left(k_{\displaystyle \scalebox{0.85}{$\hspace{0.01in}a$}}-1\right)\over 2}.$$
Since $\displaystyle \lim_{a\to +\infty}k_{\displaystyle \scalebox{0.85}{$\hspace{0.01in}a$}}=+\infty$, we can assume that $\displaystyle {k_{\displaystyle \scalebox{0.85}{$\hspace{0.01in}a$}} \left(k_{\displaystyle \scalebox{0.85}{$\hspace{0.01in}a$}}-1\right) \over 4}\geq \left(k_{\displaystyle \scalebox{0.85}{$\hspace{0.01in}a$}}+1\right)\left(q-p\right)$ and $\displaystyle {{\left( q-1\right) k_{\displaystyle \scalebox{0.85}{$\hspace{0.01in}a$}}}\over 2}\geq \scalebox{1.05}{$s_{0}$}$ for sufficiently large $a\in \mathbb{N}$. Hence we can write
$${i_{\scalebox{0.85}{$\hspace{0.01in}a$}}\over 2}\geq {n_{\displaystyle \hspace{0.01in} \scalebox{1.05}{$s_{\displaystyle \hspace{0.01in}\scalebox{0.85}{$k_{\displaystyle \scalebox{0.85}{$\hspace{0.01in}a$}}$}}$}} \over 2}\geq \scalebox{1.05}{$s_{0}$}+\left(k_{\displaystyle \scalebox{0.85}{$a$}}+1\right)\left(q-p\right)=\scalebox{1.05}{$s_{\displaystyle \hspace{0.01in}\scalebox{0.85}{$\left( k_{\displaystyle \scalebox{0.85}{$\hspace{0.01in}a$}}+1\right)$} }$}$$
and we find 
$$j_{\scalebox{0.85}{$\hspace{0.01in}a$}}= i_{\scalebox{0.85}{$\hspace{0.01in}a$}}-\scalebox{1.05}{$s_{\displaystyle \hspace{0.01in}\scalebox{0.85}{$\left( k_{\displaystyle \scalebox{0.85}{$\hspace{0.01in}a$}}+1\right)$}}$}+a>i_{\scalebox{0.85}{$\hspace{0.01in}a$}}-{i_{\scalebox{0.85}{$\hspace{0.01in}a$}}\over 2}= {i_{\scalebox{0.85}{$\hspace{0.01in}a$}}\over 2}\hspace{0.2in}\Rightarrow\hspace{0.2in} 2j_{\scalebox{0.85}{$\hspace{0.01in}a$}}>i_{\scalebox{0.85}{$\hspace{0.01in}a$}}.$$
\par Now, we will show that the inequality \ref{stability} is satisfied for a sufficiently large element of $J$. 
Let take an $a\in \mathbb{N}$ satisfiying $2j_{\scalebox{0.85}{$\hspace{0.01in}a$}}>i_{\scalebox{0.85}{$\hspace{0.01in}a$}}$. We know that $i_{\scalebox{0.85}{$\hspace{0.01in} a$}}$ is the greatest element of $m\in \mathbb{N}-I$ satisfying $ e^{\displaystyle  (c_{pq}-1)\hspace{0.0125in}{\alpha}_{\displaystyle \hspace{0.01in}\scalebox{0.85}{${n_{\displaystyle \scalebox{0.9}{$\hspace{0.01in}a$}}}$}}}\leq e^{\displaystyle c_{pq}\hspace{0.025in}\alpha_{\displaystyle\hspace{-.035in}\scalebox{0.8}{ $m$}}}$, then we can write 
$$e^{\displaystyle c_{pq}\hspace{0.025in}\alpha_{\displaystyle\hspace{-.035in}\scalebox{0.8}{ $k$}}}< e^{\displaystyle  (c_{pq}-1)\hspace{0.0125in}{\alpha}_{\displaystyle \hspace{0.01in}\scalebox{0.85}{${n_{\displaystyle \scalebox{0.9}{$\hspace{0.01in}a$}}}$}}}$$
for every $k>i_{\scalebox{0.85}{$\hspace{0.01in}a$}}$, $k\in \mathbb{N}-I$. If $2j_{\scalebox{0.85}{$\hspace{0.01in}a$}}\in \mathbb{N}-I$, then 
$$e^{\displaystyle c_{pq}\hspace{0.025in}\alpha_{\displaystyle\hspace{-.035in}\scalebox{0.85}{ $4\hspace{0.01in}j_{\scalebox{0.9}{$\hspace{0.01in}a$}}$}}}\leq e^{\displaystyle c_{pq}\hspace{0.025in}\alpha_{\displaystyle\hspace{-.035in}\scalebox{0.8}{ $2\hspace{0.01in}j_{\scalebox{0.85}{$\hspace{0.01in}a$}}$}}}<  e^{\displaystyle  (c_{pq}-1)\hspace{0.0125in}{\alpha}_{\displaystyle \hspace{0.01in}\scalebox{0.85}{${n_{\displaystyle \scalebox{0.85}{$a$}}}$}}}.$$
If $2j_{\scalebox{0.85}{$\hspace{0.01in}a$}}\in I$, then $2j_{\scalebox{0.85}{$\hspace{0.01in}a$}}+(q-p)\in \mathbb{N}-I$ and $2j_{\scalebox{0.85}{$\hspace{0.01in}a$}}+(q-p)\leq 4j_{\scalebox{0.85}{$\hspace{0.01in}a$}}$ is satisfied for a sufficiently large $a$ and we find
$$e^{\displaystyle c_{pq}\hspace{0.025in}\alpha_{\displaystyle\hspace{-.035in}\scalebox{0.8}{ $4\hspace{0.01in}j_{\scalebox{0.85}{$\hspace{0.01in}a$}}$}}}\leq e^{\displaystyle c_{pq}\hspace{0.025in}\alpha_{\displaystyle\hspace{-.035in}\scalebox{0.8}{ $2\hspace{0.01in}j_{\scalebox{0.85}{$\hspace{0.01in}a$}}$}}+(q-p)}<  e^{\displaystyle  (c_{pq}-1)\hspace{0.0125in}{\alpha}_{\displaystyle \hspace{0.01in}\scalebox{0.85}{${n_{\displaystyle \scalebox{0.85}{$a$}}}$}}}.$$
Also, we know that $i_{\scalebox{0.85}{$\hspace{0.01in}a$}}\geq j_{\scalebox{0.85}{$\hspace{0.01in}a$}}$ for every $a\in \mathbb{N}$, thus we can write
$$
d_{\displaystyle \scalebox{0.85}{$j_{\scalebox{0.85}{$a$}}$}}(U_{q}, U_{p})= e^{\displaystyle  (c_{pq}-1)\hspace{0.0125in}{\alpha}_{\displaystyle \hspace{0.01in}\scalebox{0.85}{${n_{\displaystyle \scalebox{0.85}{$\hspace{0.01in}a$}}}$}}}\leq e^{\displaystyle c_{pq}\hspace{0.025in}\alpha_{\displaystyle\hspace{-.025in}\scalebox{0.8}{ $i_{\scalebox{0.85}{$\hspace{0.01in}a$}}$}}}\leq e^{\displaystyle c_{pq}\hspace{0.025in}\alpha_{\displaystyle\hspace{0.0in}\scalebox{0.8}{$j_{\scalebox{0.85}{$\hspace{0.01in}a$}}$}}}.
$$
The above inequalites give us that
$$ e^{\mathlarger{c_{pq}\alpha_{\hspace{0.01in}4j_{\scalebox{0.85}{$a$}}}}}\leq
d_{\displaystyle \scalebox{0.85}{$j_{\scalebox{0.85}{$a$}}$}}(U_{q}, U_{p})=  e^{\displaystyle  (c_{pq}-1)\hspace{0.0125in}{\alpha}_{\displaystyle \hspace{0.01in}\scalebox{0.85}{${n_{\displaystyle \scalebox{0.85}{$a$}}}$}}}\leq e^{\displaystyle c_{pq}\hspace{0.025in}\alpha_{\displaystyle\hspace{-.04in}\scalebox{0.8}{ $j_{\scalebox{0.85}{$\hspace{0.01in}a$}}$}}}.
$$
Then, the inequality \ref{stability} is satisfied for sufficiently large element of $J$. 
\par  We now prove that the inequality \ref{stability} is satisfied for sufficiently large elements of $K$, $L$ and $M$. In order to see this, we first show that
$$n_{\displaystyle \hspace{0.01in} \scalebox{1.025}{$s_{\displaystyle \hspace{0.005in}\scalebox{0.85}{$k$}}$} }\geq 2\hspace{0.02in}\scalebox{1.05}{$s_{\displaystyle \hspace{0.005in}\scalebox{0.8}{$k$}} $}$$
for sufficently large $k\in \mathbb{N}$. We know that $n_{\displaystyle \hspace{0.01in} \scalebox{1.025}{$s_{\displaystyle \hspace{0.01in}\scalebox{0.85}{$k$}}$} }$ is on the line which has equation $x+y=q+k-2$ for every $k=0,1...$. Since the first element  of $I_{\displaystyle \scalebox{0.85}{$q+k_{\scalebox{0.85}{$\hspace{0.01in}a$}}-2$}}$ is less than $n_{\displaystyle \hspace{0.01in} \scalebox{1.025}{$s_{\displaystyle \hspace{0.01in}\scalebox{0.85}{$k$}}$} }$, then we can write
$$n_{\displaystyle \hspace{0.01in} \scalebox{1.025}{$s_{\displaystyle \hspace{0.005in}\scalebox{0.85}{$k$}}$}}\geq {\left(q+k-2\right)\left(q+k-1\right)\over 2}= {{\left( q-2\right) \left(q-1\right)}\over 2}+\left( q-1\right) k+{k\left(k-1\right)\over 2}.$$
The inequalities 
$${k.(k-1)\over 4} \geq k(q-p) \hspace{0.25in} \textnormal{and} \hspace{0.25in} (q-1)k\geq 2(\scalebox{1.05}{$s_{0}$}+1)$$
hold for a sufficiently large $k$. Then we find 
$$n_{\displaystyle \hspace{0.01in} \scalebox{1.025}{$s_{\displaystyle \hspace{0.005in}\scalebox{0.85}{$k$}}$} }\geq 2\hspace{0.01in}(\scalebox{1.05}{$s_{0}$}+k(q-p)+1)= 2\hspace{0.02in}\scalebox{1.025}{$s_{\displaystyle \hspace{0.005in}\scalebox{0.85}{$k$}}$}$$
for a sufficiently large $k\in \mathbb{N}$.
\par Now we show that the inequality \ref{stability} is satisfied for sufficiently large element of $\displaystyle K$.  Let take an $n\in K$. Then, there exist a $a\in \mathbb{N}$ and a $k\in \mathbb{N}$ satisfying $k_{\displaystyle 	\hspace{0.01in}\scalebox{0.85}{$a$}}+1\leq k\leq k_{\displaystyle \hspace{0.01in}\scalebox{0.85}{$(a+1)$}}-1$ such that
$n_{\displaystyle \hspace{0.01in} \scalebox{1.025}{$s_{\displaystyle \hspace{0.01in}\scalebox{0.85}{$k$}}$} }-\scalebox{1.05}{$s_{\displaystyle \hspace{0.01in}\scalebox{0.85}{$k$}}$}+a \leq n\leq n_{\displaystyle \hspace{0.01in} \scalebox{1.025}{$s_{\displaystyle \hspace{0.01in}\scalebox{0.8}{$(k+1)$}} $}}-\scalebox{1.05}{$s_{\displaystyle \hspace{0.01in}\scalebox{0.8}{$(k+1)$}}$}+a-1
$
and 
$$d_{n}\left(U_{q},U_{p}\right)= e^{\displaystyle {c_{pq}\hspace{0.015in}\alpha_{\displaystyle \scalebox{0.9}{$\hspace{0.005in}n+\scalebox{1.05}{$s_{\displaystyle \scalebox{0.8}{$(k+1)$}}$}-a$}}}}.$$
Since $n_{\displaystyle \hspace{0.01in} \scalebox{1.025}{$s_{\displaystyle \hspace{0.01in}\scalebox{0.85}{$k$}}$} }\geq 2\hspace{0.01in}\scalebox{1.05}{$s_{\displaystyle \hspace{0.015in}\scalebox{0.8}{$k$}}$}$ for a sufficiently large $k\in \mathbb{N}$ and $\scalebox{1.05}{$s_{\displaystyle \hspace{0.01in}\scalebox{0.8}{$(k+1)$}}$}-\scalebox{1.05}{$s_{\displaystyle \hspace{0.01in}\scalebox{0.8}{$k$}}$}=q-p$ for all $k\in \mathbb{N}$, we can write
$$\scalebox{1.025}{$s_{\displaystyle \hspace{0.01in}\scalebox{0.85}{$k$}}$}\leq n_{\displaystyle \hspace{0.01in} \scalebox{1.025}{$s_{\displaystyle \hspace{0.01in}\scalebox{0.8}{$k$}}$}}-\scalebox{1.025}{$s_{\displaystyle \hspace{0.015in}\scalebox{0.8}{$k$}}$}+a\leq n \hspace{0.2in} \Rightarrow \hspace{0.2in}n+\scalebox{1.025}{$s_{\displaystyle \hspace{0.015in}\scalebox{0.8}{$(k+1)$}}$}-a\leq 2n. $$
for sufficiently large $a$. Then, we obtain 
$$d_{n}\left(U_{q},U_{p}\right)= e^{\displaystyle {c_{pq}\hspace{0.015in}\alpha_{\displaystyle \scalebox{0.9}{$\hspace{0.01in}n+\scalebox{1.025}{$s_{\displaystyle \scalebox{0.8}{$(k+1)$}}$}-a$}}}}\geq e^{\mathlarger{c_{pq}\hspace{0.015in}\alpha_{\displaystyle \scalebox{0.8}{$\hspace{0.01in}2n$}}}}\geq e^{\mathlarger{c_{pq}\hspace{0.015in}\alpha_{\displaystyle \scalebox{0.8}{$\hspace{0.01in}4n$}}}}$$
and always we have
$$d_{n}\left(U_{q},U_{p}\right)= e^{\displaystyle {c_{pq}\hspace{0.015in}\alpha_{\displaystyle \scalebox{0.9}{$\hspace{0.01in}n+\scalebox{1.025}{$s_{\displaystyle \scalebox{0.8}{$(k+1)$}}$}-a$}}}}\leq e^{\mathlarger{c_{pq}\hspace{0.015in}\alpha_{\displaystyle \scalebox{0.8}{$\hspace{0.01in}n$}}}}$$
since $\alpha$ is increasing. Therefore, the inequality \ref{stability} is satisfied for sufficiently large elemets of $K$.
\par Now, we will show that the inequality \ref{stability} is satisfied for a sufficiently large  element of $\displaystyle L$.  Let us take a $n\in L$. Then, there is an $a\in \mathbb{N}$ such that
$$j_{\scalebox{0.85}{$\hspace{0.01in}a$}}+1\leq n\leq  n_{\displaystyle\hspace{0.01in}{\scalebox{1.025}{$s_{\displaystyle \scalebox{0.8}{$(k_{\scalebox{0.8}{$a$}}+1)$}}$}}}-\scalebox{1.025}{$s_{\displaystyle \hspace{0.01in} \scalebox{0.8}{$(k_{\scalebox{0.85}{$\hspace{0.01in}a$}}+1)$}}$}+a-1$$
and 
$$d_{n}\left(U_{q},U_{p}\right)= e^{\displaystyle {c_{pq}\hspace{0.015in}\alpha_{\displaystyle \scalebox{0.9}{$\hspace{0.01in}n+\scalebox{1.025}{$s_{\displaystyle \scalebox{0.9}{$\hspace{0.01in}(k_{\scalebox{0.85}{$\hspace{0.01in}a$}}+1)$}}$}-a$}}}}.$$
Since $\scalebox{1.025}{$s_{\displaystyle \hspace{0.01in}\scalebox{0.8}{$k_{\scalebox{0.85}{$\hspace{0.01in}a$}}$}}$} \leq n_{\displaystyle \hspace{0.01in}{\scalebox{1.025}{$s_{\displaystyle \hspace{0.01in} \scalebox{0.8}{$k_{\scalebox{0.85}{$\hspace{0.01in}a$}}$}}$}}}-\scalebox{1.05}{$s_{\displaystyle \hspace{0.01in}\scalebox{0.8}{$k_{\scalebox{0.85}{$\hspace{0.01in}a$}}$}}$}+a\leq j_{\scalebox{0.85}{$\hspace{0.01in}a$}}+1\leq n$
and $ n+\scalebox{1.05}{$s_{\displaystyle \hspace{0.01in} \scalebox{0.8}{$(k_{\scalebox{0.85}{$\hspace{0.01in}a$}}+1)$}}$}-a\leq 2n$ for a sufficiently large $n$, then we find
$$d_{n}\left(U_{q},U_{p}\right)= e^{\displaystyle {c_{pq}\hspace{0.015in}\alpha_{\displaystyle \scalebox{0.9}{$\hspace{0.01in}n+\scalebox{1.025}{$s_{\displaystyle \scalebox{0.9}{$\hspace{0.01in}(k_{\scalebox{0.85}{$\hspace{0.01in}a$}}+1)$}}$}-a$}}}}\geq e^{\displaystyle c_{pq}\hspace{0.025in}\alpha_{\displaystyle\hspace{-.025in}\scalebox{0.8}{ $\hspace{0.01in}2n$}}}\geq  e^{\displaystyle c_{pq}\hspace{0.025in}\alpha_{\displaystyle\hspace{-.025in}\scalebox{0.8}{ $\hspace{0.01in}4n$}}},$$
and always we have 
$$d_{n}\left(U_{q},U_{p}\right)= e^{\displaystyle {c_{pq}\hspace{0.015in}\alpha_{\displaystyle \scalebox{0.9}{$\hspace{0.01in}n+\scalebox{1.025}{$s_{\displaystyle \scalebox{0.9}{$\hspace{0.01in}(k_{\scalebox{0.85}{$\hspace{0.01in}a$}}+1)$}}$}-a$}}}}\leq e^{\mathlarger{c_{pq}\hspace{0.015in}\alpha_{\displaystyle \scalebox{0.8}{$\hspace{0.01in}n$}}}}$$
since $\alpha$ is increasing. Therefore, the inequality \ref{stability} is satisfied for sufficiently large element of $L$.
\par Now we will show that the inequality \ref{stability} is satisfied for a sufficiently large element of $\displaystyle M$. If $n\in M$, then there is an $a\in \mathbb{N}$
$$n_{\displaystyle \hspace{0.01in}{\scalebox{1.025}{$s_{\displaystyle \hspace{0.01in} \scalebox{0.9}{$k_{\scalebox{0.85}{$\hspace{0.01in}a$}}$}}$}}}- \scalebox{1.025}{$s_{\displaystyle \hspace{0.01in}\scalebox{0.9}{$(k_{\scalebox{0.85}{$\hspace{0.01in}a$}}+1)$}}$}+a\leq n\leq  j_{\scalebox{0.85}{$\hspace{0.01in}a$}}-1$$
and
$$d_{n}\left(U_{q},U_{p}\right)= e^{\displaystyle {c_{pq}\hspace{0.015in}\alpha_{\displaystyle \scalebox{0.9}{$\hspace{0.01in}n+\scalebox{1.1}{$s_{\displaystyle \scalebox{0.85}{$\hspace{0.01in}(k_{\scalebox{0.85}{$\hspace{0.01in}a$}}+1)$}}$}-(a-1)$}}}}.$$
Again we can write $\scalebox{1.1}{$s_{\displaystyle \hspace{0.01in} \scalebox{0.8}{$k_{\scalebox{0.85}{$\hspace{0.01in}a$}}$}}$} \leq n_{\displaystyle \scalebox{1.05}{$s_{\displaystyle \hspace{0.01in} \scalebox{0.8}{$k_{\scalebox{0.85}{$\hspace{0.01in}a$}}$}}$}}-\scalebox{1.1}{$s_{\displaystyle \hspace{0.01in} \scalebox{0.8}{$k_{\scalebox{0.85}{$\hspace{0.01in}a$}}$}}$}+a\leq n$ and $ n_{\displaystyle \scalebox{1.05}{$s_{\displaystyle \hspace{0.01in} \scalebox{0.8}{$k_{\scalebox{0.85}{$\hspace{0.01in}a$}}$}}$}}+\scalebox{1.05}{$s_{\displaystyle \hspace{0.01in} \scalebox{0.8}{$k_{\scalebox{0.85}{$\hspace{0.01in}a$}}$}}$}-a+1\leq 2n$ for a sufficiently large $a$. Hence we find
$$d_{n}\left(U_{q},U_{p}\right)=e^{\displaystyle {c_{pq}\hspace{0.015in}\alpha_{\displaystyle \scalebox{0.9}{$\hspace{0.01in} n+\scalebox{1.05}{$s_{\displaystyle \hspace{0.01in}  \scalebox{0.8}{$(k_{\scalebox{0.85}{$\hspace{0.01in}a$}}+1)$}}$}-(a-1)$}}}}\geq	e^{\mathlarger{c_{pq}\hspace{0.015in}\alpha_{\displaystyle \scalebox{0.8}{$\hspace{0.01in}2n$}}}}\geq 	e^{\mathlarger{c_{pq}\hspace{0.015in}\alpha_{\displaystyle \scalebox{0.8}{$\hspace{0.01in}4n$}}}}$$
and always we have 
$$d_{n}\left(U_{q},U_{p}\right)=e^{\displaystyle {c_{pq}\hspace{0.015in}\alpha_{\displaystyle \scalebox{0.9}{$\hspace{0.01in} n+\scalebox{1.05}{$s_{\displaystyle \hspace{0.01in}  \scalebox{0.8}{$(k_{\scalebox{0.85}{$\hspace{0.01in}a$}}+1)$}}$}-(a-1)$}}}} \leq e^{\mathlarger{c_{pq}\hspace{0.015in}\alpha_{\displaystyle \scalebox{0.8}{$\hspace{0.01in}n$}}}}$$
since $\alpha$ is increasing. Therefore, the inequality \ref{stability} is satisfied for a sufficiently large element of $M$. This completes the proof.
\end{proof}
~\\$\vspace{-0.35in}$~\\
\subsection{The diametral dimension and the approximate diametral dimension of an element of the family $\bf \mathcal{K}$ parameterized by a sequence $\boldsymbol{\alpha}$} ~\\
\par As a consequence of Theorem \ref{stabilitycase}, we  will compute the diametral dimension and the approximate diametral dimension of an element $\bf \mathcal{K}_{\alpha}$ of the family $\bf \mathcal{K}$ which is parameterized by a stable sequence $\alpha$.
\begin{cor}\label{sc}  Let $\bf \mathcal{K}_{\alpha}$ be an element of the family $\bf \mathcal{K}$ which is parameterized by a stable sequence $\alpha$. Then, $\displaystyle \Delta(  \mathcal{K}_{\alpha} )=\Delta(\Lambda_{1}\left(\alpha_{n}\right))$ and $\displaystyle\delta(  \mathcal{K}_{\alpha} )=\delta(\Lambda_{1}\left(\alpha_{n}\right))$.
\end{cor}
\begin{proof} From Theorem \ref{stabilitycase}, we have
$$\Delta(\Lambda_{1}(\alpha_{ \displaystyle \scalebox{0.8}{$n$}}))\subseteq\Delta ( \mathcal{K}_{\alpha})\subseteq \Delta(\Lambda_{1}(\alpha_{\displaystyle \hspace{0.01in} \scalebox{0.8}{$4n$}}))$$
and 
$$\delta(\Lambda_{1}(\alpha_{\displaystyle \scalebox{0.8}{$4n$}}))\subseteq\delta ( \mathcal{K}_{\alpha})\subseteq \delta(\Lambda_{1}(\alpha_{\displaystyle \scalebox{0.8}{$n$}})).$$
On the other hand, $\Lambda_{1}(\alpha_{\displaystyle \scalebox{0.8}{$\hspace{0.01in}n$}})\cong \Lambda_{1}(\alpha_{\displaystyle \scalebox{0.8}{$\hspace{0.01in} 4n$}})$ since $\alpha$ is stable. Then 
$\Delta( \mathcal{K}_{\alpha})=\Delta(\Lambda_{1}\left(\alpha_{\displaystyle \scalebox{0.8}{$n$}}\right))$ and $\delta( \mathcal{K}_{\alpha} )=\delta(\Lambda_{1}\left(\alpha_{\displaystyle \scalebox{0.8}{$n$}}\right))$.
\end{proof}
\par Now we will prove that  $\Delta( \mathcal{K}_{\alpha} )=\Delta\left(\Lambda_{1}\left(\alpha_{\displaystyle \scalebox{0.8}{$n+1$}}\right)\right)$ and $\delta( \mathcal{K}_{\alpha})\neq \delta\left(\Lambda_{1}\left(\alpha_{\displaystyle \scalebox{0.8}{$n+1$}}\right)\right)$ for an element $\bf \mathcal{K}_{\alpha}$ of the family $\bf \mathcal{K}$ which is parameterized by an unstable sequence $\alpha$. Besides, we will show that all regular elements of the family $\bf \mathcal{K}$ are parameterized by an unstable sequence $\alpha$.
\begin{prop}\label{EDD}  Let $\bf \mathcal{K}_{\alpha}$ be an element of the family $\bf \mathcal{K}$ which is parameterized by an unstable sequence $\alpha$. Then, $\displaystyle \Delta( \mathcal{K}_{\alpha} )=\Delta(\Lambda_{1}\left(\alpha_{\displaystyle \scalebox{0.8}{$n+1$}}\right))$.\end{prop}
\begin{proof} We can calculate Kolmogorov diameters as in the previous determined for every $p$ and $q>p$. Since $\alpha $ is unstable, then there exists an $a_{0}\in \mathbb{N}$ such that for all $a\geq a_{0}$, there is no $m>n_{\displaystyle \scalebox{0.8}{$a$}}$, $m\in \mathbb{N}$ satisfying
$\alpha_{\displaystyle\hspace{-.025in}\scalebox{0.8}{ $m$}}\leq A_{\displaystyle \scalebox{0.8}{$pq$}}\hspace{0.025in}\alpha_{\hspace{0.01in}\scalebox{0.85}{${n_{\displaystyle \scalebox{0.85}{$a$}}}$}}.$
Now, we examine closely the indices replaced the term $ \displaystyle { e^{\displaystyle  (c_{pq}-1)\hspace{0.0125in}{\alpha}_{\displaystyle \hspace{0.01in}\scalebox{0.9}{${n_{\displaystyle \scalebox{0.9}{$\hspace{0.01in}a_{\displaystyle \scalebox{0.8}{$\hspace{0.01in}0$}}$}}}$}}}}$.
We know that
$$d_{\displaystyle j_{\displaystyle \scalebox{0.85}{$\hspace{0.01in}(a_{0}-1)$}}} (U_{q}, U_{p})={ e^{\displaystyle  (c_{pq}-1)\hspace{0.0125in}{\alpha}_{\displaystyle \hspace{0.01in}\scalebox{0.9}{${{n_{\displaystyle \scalebox{0.85}{$\hspace{0.01in}(a_{\displaystyle \scalebox{0.8}{$\hspace{0.01in}0$}}-1)$}}}}$}}}}$$
where $j_{\displaystyle 
\scalebox{0.85}{$\hspace{0.01in}(a_{0}-1)$}}=i_{\displaystyle 
\scalebox{0.85}{$\hspace{0.01in}(a_{0}-1)$}}-\scalebox{1.1}{$s_{\displaystyle \hspace{0.01in}\scalebox{0.8}{$(a_{0}-1)$}}$}+\displaystyle {a_{0}-2}$.
Since $\alpha_{\displaystyle \hspace{0.01in} \scalebox{0.95}{$i_{\displaystyle \hspace{0.01in}\scalebox{0.85}{$\hspace{0.01in}(a_{\hspace{0.01in}0}-1)$}} $}}\leq A_{\displaystyle \scalebox{0.8}{$pq$}}\hspace{0.025in} {\alpha}_{\displaystyle \hspace{0.01in}\scalebox{0.85}{${{n_{\displaystyle \scalebox{0.95}{$\hspace{0.01in}(a_{\displaystyle \scalebox{0.8}{$0$}}-1)$}}}}$}}$ and there is no $m >n_{\hspace{0.01in} \displaystyle{\scalebox{0.9}{$\hspace{0.01in}a_{0}$}}}$ satisfying $\alpha_{\displaystyle\hspace{-.025in}\scalebox{0.8}{ $m$}}\leq A_{\displaystyle \scalebox{0.8}{$pq$}}\hspace{0.025in}\alpha_{\hspace{0.01in}\scalebox{0.9}{${n_{\displaystyle \scalebox{0.9}{$a_{0}$}}}$}}$, then we find $i_{\displaystyle 
\scalebox{0.85}{$(a_{0}-1)$}}< n_{\hspace{0.01in} \displaystyle{\scalebox{0.9}{$a_{0}$}}}$. This gives that for all ${j_{\displaystyle \scalebox{0.85}{$(a_{\displaystyle \scalebox{0.8}{$0$}}-1)$}}} \leq n\leq {n_{\displaystyle \scalebox{0.85}{$a_{\displaystyle \scalebox{0.75}{$\hspace{0.01in}0$}}$}}}-2$, 
$$d_{n}(U_{q}, U_{p})= e^{\displaystyle c_{pq}\hspace{0.025in}\alpha_{\displaystyle \scalebox{0.9}{$n+1$}}}.$$
Besides, we obtain that the sequence $\displaystyle \left(a_{p,n}\over a_{q,n}\right)_{n\in \mathbb{N}}$ has decreasing order starting from the indices $j_{\displaystyle \scalebox{0.85}{$(a_{\displaystyle \scalebox{0.8}{$0$}}-1)$}} +1$, since for every $a\geq a_{0}$, there is no $n>n_{\hspace{0.01in} \displaystyle a_{0}}$ satisfying $\alpha_{\displaystyle \scalebox{0.9}{$n$}}\leq A_{pq}\hspace{0.025in} \alpha_{\hspace{0.01in} \displaystyle \scalebox{0.9}{$n_{\hspace{0.01in} \displaystyle a}$}}$. Then, we have for all $a\geq a_{0}$ 
$$d_{\hspace{0.01in} \displaystyle \scalebox{0.9}{$n_{\hspace{0.01in} \displaystyle \scalebox{0.9}{$a$}-1}$}}(U_{q}, U_{p})=e^{\mathlarger{(c_{pq}-1)\alpha_{\hspace{0.01in} \displaystyle\scalebox{0.9}{$n_{\displaystyle \scalebox{0.8}{$a$}}$}}}}$$
and for all $m\geq j_{\displaystyle \scalebox{0.9}{$(a_{\displaystyle \scalebox{0.8}{$0$}}-1)$}}$, $m\in \mathbb{N}-I$
$$d_{m}(U_{q}, U_{p})=e^{\displaystyle c_{pq}\hspace{0.025in}\alpha_{\displaystyle\hspace{-.025in}\scalebox{0.8}{ $m+1$}}}.$$ 
Since $\displaystyle d_{n}(U_{q}, U_{p})\leq e^{\displaystyle c_{pq}\hspace{0.025in}\alpha_{\displaystyle\hspace{-.025in}\scalebox{0.8}{ $n+1$}}}$ for every $n\in \mathbb{N}$, then we find $\Delta( \mathcal{K}_{\alpha})\supseteq\Delta(\Lambda_{1}(\alpha_{\displaystyle \scalebox{0.8}{$n+1$}}))$. 
\par For the other direction, let us take a sequence $\left( x_{n}\right)_{n\in \mathbb{N}}\in  \Delta(\bf \mathcal{K}_{\alpha})$, an $\varepsilon>0$ and a $p\in \mathbb{N}$ satisfying $\displaystyle {1\over p}<\varepsilon$. We will show that
$$\sup_{n\in \mathbb{N}}\left|x_{n}\right|e^{\displaystyle -\varepsilon\hspace{0.025in}\alpha_{\displaystyle\hspace{-.025in}\scalebox{0.8}{ $n+1$}}}<+\infty.$$
Since $\left( x_{n}\right)_{n\in \mathbb{N}}\in  \Delta(\bf \mathcal{K}_{\alpha})$, there exist a $q>p$ and $M_{1}>0$ satisfying 
$$\sup_{n\in \mathbb{N}}\left|x_{n}\right|d_{n}\left(U_{p}, U_{q}\right)<M_{1}.$$
Let us define $\displaystyle I=\bigcup_{\displaystyle \scalebox{0.8}{$p\leq s<q$}} I_{\displaystyle \scalebox{0.9}{$s$}}$. For sufficiently large  $n\in \mathbb{N}-I$, we can write
$$\left|x_{n}\right|e^{\displaystyle {-\varepsilon  \hspace{0.01in}\alpha_{\displaystyle \scalebox{0.8}{$n+1$}}}}\leq \left|x_{n}\right|d_{n}\left(U_{q},U_{p}\right)=e^{\displaystyle {c_{pq} \hspace{0.01in}\alpha_{\displaystyle \scalebox{0.8}{$n+1$}}}}\leq M_{1}$$
since $c_{pq}\geq -\varepsilon$. Therefore, the sequence $\left|x_{n}\right|e^{\displaystyle {-\varepsilon  \hspace{0.01in}\alpha_{\displaystyle \scalebox{0.8}{$n+1$}}}}$ is bounded on the set $\mathbb{N}-I$. If we show that  $\left|x_{n}\right|e^{\displaystyle {-\varepsilon  \hspace{0.01in}\alpha_{\displaystyle \scalebox{0.8}{$n+1$}}}}$ is also bounded on $I$, then we will find that $\left(x_{n}\right)_{n\in \mathbb{N}}\in \Delta(\Lambda_{1}\left(\alpha_{\displaystyle \scalebox{0.8}{$n+1$}}\right))$. Let take another $p_{0}>q$, then there exist a $q_{0}$ and $M_{2}>0$ such that  
$$\sup_{n\in \mathbb{N}}\left|x_{n}\right|d_{n}\left(U_{\displaystyle \scalebox{0.8}{$q_{0}$}}, U_{\displaystyle \scalebox{0.8}{$p_{0}$}}\right)<M_{2}.$$
Let us define $\displaystyle J =\bigcup_{\displaystyle\scalebox{0.8}{$p_{0}\leq s<q_{0}$}} I_{\displaystyle \scalebox{0.8}{$s$}}$. Since $c_{\displaystyle \scalebox{0.8}{$p_{0}$},\displaystyle \scalebox{0.8}{$q_{0}$}}\geq -\varepsilon$, we find
$$\left|x_{n}\right|e^{\displaystyle {-\varepsilon  \hspace{0.01in}\alpha_{\displaystyle \scalebox{0.8}{$n+1$}}}}\leq \left|x_{n}\right|d_{n}\left(U_{\displaystyle \scalebox{0.8}{$q_{0}$}},U_{\displaystyle \scalebox{0.8}{$p_{0}$}}\right)=e^{\displaystyle {c_{\displaystyle \scalebox{0.8}{$p_{0}$},\displaystyle \scalebox{0.8}{$q_{0}$}} \hspace{0.01in}\alpha_{\displaystyle \scalebox{0.8}{$n+1$}}}}\leq M_{2}$$
for sufficently large $n+1\in \mathbb{N}-J$. Also, it is easy to see that $I\subset \mathbb{N}-J$. Then, the above inequalities give us that 
$$\left|x_{n}\right|e^{\displaystyle {-\varepsilon  \hspace{0.01in}\alpha_{\displaystyle \scalebox{0.8}{$n+1$}}}}\leq  M_{2}.$$
for all $n\in I$. Hence, the sequence  $\left|x_{n}\right|e^{\displaystyle {-\varepsilon  \hspace{0.01in}\alpha_{\displaystyle \scalebox{0.8}{$n+1$}}}}$ is also bounded on $I$. Therefore, we find  
$$\sup_{n\in \mathbb{N}}\left|x_{n}\right|e^{\displaystyle {-\varepsilon  \hspace{0.01in}\alpha_{\displaystyle \scalebox{0.8}{$n+1$}}}}<+\infty$$ 
and  $\left(x_{n}\right)_{n\in \mathbb{N}}\in \Delta(\Lambda_{1}\left(\alpha_{\displaystyle \scalebox{0.8}{$n+1$}}\right))$. This says that $\Delta( \mathcal{K}_{\alpha} )=\Delta(\Lambda_{1}\left(\alpha_{\displaystyle \scalebox{0.8}{$n+1$}}\right))$.
\end{proof}	
\begin{prop}\label{EADD}  Let $\bf \mathcal{K}_{\alpha}$ be an element of the family $\bf \mathcal{K}$ which is parameterized by an unstable sequence $\alpha$. Then, $\displaystyle\delta(  \mathcal{K}_{\alpha})\neq \delta(\Lambda_{1}\left(\alpha_{\displaystyle \scalebox{0.8}{$n+1$}}\right))$.
\end{prop}
\begin{proof} In the proof of the previous proposition, we show that if $\alpha$ is unstable, then for all $p\in \mathbb{N}$ and $q>p$, there is a $a_{0}\in \mathbb{N}$ such that for all $a\geq a_{0}$
$$d_{\hspace{0.01in} \displaystyle \scalebox{0.9}{$n_{\hspace{0.01in} \displaystyle \scalebox{0.9}{$a$}-1}$}}(U_{q}, U_{p})=e^{\mathlarger{(c_{pq}-1)\alpha_{\hspace{0.01in} \displaystyle\scalebox{0.9}{$n_{\displaystyle \scalebox{0.8}{$a$}}$}}}},$$
so the last equality holds except for finitely many numbers of elements of $I$.
Then we have
$${\varepsilon_{\displaystyle \scalebox{0.9}{$n_{\displaystyle \scalebox{0.8}{$a$}}-1$}} \left(p,q\right)\over \alpha_{\displaystyle \scalebox{0.9}{$n_{\displaystyle \scalebox{0.8}{$a$}}$}}}= 1-c_{p,q}$$
and
$$\limsup_{\displaystyle \scalebox{0.8}{$a\in \mathbb{N}$}}{\varepsilon_{\displaystyle \scalebox{0.9}{$n_{\displaystyle \scalebox{0.8}{$a$}}-1$}} \left(p,q\right)\over \alpha_{\displaystyle \scalebox{0.9}{$n_{\displaystyle \scalebox{0.8}{$a$}}$}}}= 1-c_{p,q}\hspace{0.25in}\Rightarrow \hspace{0.25in}\inf_{p}\sup_{q} \limsup_{n\in \mathbb{N}}{\varepsilon_{\displaystyle \scalebox{0.85}{$n$}}\left(p,q\right)\over \alpha_{\displaystyle \scalebox{0.85}{$n+1$}}}>0.$$
By Proposition \ref{AA}, we have $\delta\left(\bf \mathcal{K}_{\alpha} \right)\neq \delta\left(\Lambda_{1}\left(\alpha_{\displaystyle \scalebox{0.85}{$n+1$}}\right)\right)$. 
\end{proof} 
\begin{remark} Proposition \ref{EDD} and Proposition \ref{EADD} shows that Question \ref{qi2} has a negative answer for the elements of the family $\bf \mathcal{K}$ which is parametrized by an unstable exponent sequence.
\end{remark}
\par Now, we will show that all regular elements of the family $\bf \mathcal{K}$ are parameterized by an unstable sequence $\alpha$. 
\par Let $\bf \mathcal{K}_{\alpha}$ be an element of the family $\bf \mathcal{K}$ parameterized by an exponent sequence $\alpha$ and $n\in I_{\displaystyle \scalebox{0.85}{$s$}}$, $s\in \mathbb{N}$. Then, there exist two cases for $n+1$: $n+1\in I_{\displaystyle \scalebox{0.85}{$s+1$}}$ or $n+1\in I_{1}$.
~\\$\vspace{-.2in}$~\\ $\vspace{-.1in}$~\\ We assume ${n+1\in I_{\displaystyle \scalebox{0.85}{$s+1$}}}$: For this case, $\displaystyle n+1\geq {(s+1)(s+2)\over 2}\geq s+1.$
\begin{itemize}
\item[\textbf{i)}] For ${k+1\leq s}$, we have 	$a_{\displaystyle\hspace{0.01in} \scalebox{0.8}{$k,n$}}= e^{\displaystyle{-\scalebox{1.25}{${1\over k }$}\hspace{0.01in}{\displaystyle \alpha_{\displaystyle \scalebox{0.85}{$n$}}}}},$
$a_{\displaystyle\hspace{0.01in} \scalebox{0.8}{$k+1,n$}}= e^{\displaystyle{-\scalebox{1.25}{${1\over k+1 }$}\hspace{0.01in}{\displaystyle \alpha_{\displaystyle \scalebox{0.85}{$n$}}}}},$ 	
$a_{\displaystyle\hspace{0.01in} \scalebox{0.8}{$k,n+1$}}= e^{\displaystyle{-\scalebox{1.25}{${1\over k }$}\hspace{0.01in}{\displaystyle \alpha_{\displaystyle \scalebox{0.85}{$n+1$}}}}}$
$a_{\displaystyle\hspace{0.01in} \scalebox{0.8}{$k+1,n+1$}}= e^{\displaystyle{-\scalebox{1.25}{${1\over k+1 }$}\hspace{0.01in}{\displaystyle \alpha_{\displaystyle \scalebox{0.85}{$n+1$}}}}}$. Since $\alpha$ is increasing, the inequality 
$$\displaystyle { a_{\displaystyle\hspace{0.01in} \scalebox{0.8}{$k+1,n$}}\over  a_{\displaystyle\hspace{0.01in} \scalebox{0.8}{$k,n$}}}=e^{\displaystyle {\scalebox{0.8}{$\displaystyle \left({1\over k}-{1\over k+1}\right)$}}\hspace{0.01in}{\displaystyle\alpha_{\displaystyle \scalebox{0.8}{$n$}}}}
\leq e^{\displaystyle {\scalebox{0.8}{$\displaystyle \left({1\over k}-{1\over k+1}\right)$}}\hspace{0.01in}{\displaystyle\alpha_{\displaystyle \scalebox{0.8}{$n+1$}}}}={ a_{\displaystyle\hspace{0.01in} \scalebox{0.8}{$k+1,n+1$}}\over  a_{\displaystyle\hspace{0.01in} \scalebox{0.8}{$k,n+1$}}}$$
holds in this case. 
\item[\textbf{ii)}] For ${k\geq s+1}$, we have
$a_{\displaystyle\hspace{0.01in} \scalebox{0.8}{$k,n$}}= e^{{\scalebox{0.8}{$ \displaystyle\left(-{1\over k }+1\right)$}\hspace{0.01in}{\displaystyle \alpha_{\displaystyle \scalebox{0.85}{$n$}}}}},$  
$a_{\displaystyle\hspace{0.01in} \scalebox{0.8}{$k+1,n$}}= e^{{\scalebox{0.8}{$\displaystyle\left(-{1\over k+1 }+1\right)$}\hspace{0.01in}{\displaystyle \alpha_{\displaystyle \scalebox{0.85}{$n$}}}}},$
$a_{\displaystyle\hspace{0.01in} \scalebox{0.8}{$k,n+1$}}= e^{{\scalebox{0.8}{$ \displaystyle \left(-{1\over k }+1\right)$}\hspace{0.01in}{\displaystyle \alpha_{\displaystyle \scalebox{0.85}{$n+1$}}}}}$, 
$a_{\displaystyle\hspace{0.01in} \scalebox{0.8}{$k+1,n+1$}}= e^{{\scalebox{0.8}{$\displaystyle-{1\over k+1 }$}\hspace{0.01in}{\displaystyle \alpha_{\displaystyle \scalebox{0.85}{$n+1$}}}}}$. Since $\alpha$ is increasing, the inequality
$$\displaystyle { a_{\displaystyle\hspace{0.01in} \scalebox{0.8}{$k+1,n$}}\over  a_{\displaystyle\hspace{0.01in} \scalebox{0.8}{$k,n$}}}=e^{\displaystyle {\scalebox{0.8}{$\displaystyle \left({1\over k}-{1\over k+1}\right)$}}\hspace{0.01in}{\displaystyle\alpha_{\displaystyle \scalebox{0.8}{$n$}}}}
\leq e^{\displaystyle {\scalebox{0.8}{$\displaystyle \left({1\over k}-{1\over k+1}\right)$}}\hspace{0.01in}{\displaystyle\alpha_{\displaystyle \scalebox{0.8}{$n+1$}}}}={ a_{\displaystyle\hspace{0.01in} \scalebox{0.8}{$k+1,n+1$}}\over  a_{\displaystyle\hspace{0.01in} \scalebox{0.8}{$k,n+1$}}}$$
holds in this case.
\item[\textbf{iii)}] For ${k=s}$, we have
$a_{\displaystyle\hspace{0.01in} \scalebox{0.8}{$k,n$}}= e^{\displaystyle{-\scalebox{1.25}{${1\over k }$}\hspace{0.01in}{\displaystyle \alpha_{\displaystyle \scalebox{0.85}{$n$}}}}},$
$a_{\displaystyle\hspace{0.01in} \scalebox{0.8}{$k+1,n$}}= e^{{\scalebox{0.8}{$\displaystyle \left(-{1\over k+1 }+1\right)$}\hspace{0.01in}{\displaystyle \alpha_{\displaystyle \scalebox{0.85}{$n$}}}}},$
$a_{\displaystyle\hspace{0.01in} \scalebox{0.8}{$k,n+1$}}= e^{\displaystyle{-\scalebox{1.25}{${1\over k }$}\hspace{0.01in}{\displaystyle \alpha_{\displaystyle \scalebox{0.85}{$n+1$}}}}}$,
$a_{\displaystyle\hspace{0.01in} \scalebox{0.8}{$k+1,n+1$}}= e^{\displaystyle{-\scalebox{1.25}{${1\over k+1 }$}\hspace{0.01in}{\displaystyle \alpha_{\displaystyle \scalebox{0.85}{$n+1$}}}}}.$
Then, these give that
~\\~\\$\displaystyle { a_{\displaystyle\hspace{0.01in} \scalebox{0.8}{$k+1,n$}}\over  a_{\displaystyle\hspace{0.01in} \scalebox{0.8}{$k,n$}}}=e^{\displaystyle {\scalebox{0.8}{$\displaystyle \left({1\over k}-{1\over k+1}+1\right)$}}\hspace{0.01in}{\displaystyle\alpha_{\displaystyle \scalebox{0.8}{$n$}}}}			
\hspace{0.25in}\textnormal{and}\hspace{0.25in} {a_{\displaystyle\hspace{0.01in} \scalebox{0.8}{$k+1,n+1$}}\over  a_{\displaystyle\hspace{0.01in} \scalebox{0.8}{$k,n+1$}}}=e^{\displaystyle {\scalebox{0.8}{$\displaystyle \left({1\over k}-{1\over k+1}\right)$}}\hspace{0.01in}{\displaystyle\alpha_{\displaystyle \scalebox{0.8}{$n+1$}}}.}$~\\~\\
In this case, the regularity condition $ \displaystyle { a_{\displaystyle\hspace{0.01in} \scalebox{0.8}{$k+1,n$}}\over  a_{\displaystyle\hspace{0.01in} \scalebox{0.8}{$k,n$}}}\leq  { a_{\displaystyle\hspace{0.01in} \scalebox{0.8}{$k+1,n+1$}}\over  a_{\displaystyle\hspace{0.01in} \scalebox{0.8}{$k,n+1$}}}$ is equivalent to the following inequality:
\begin{equation}\label{regular}
\hspace{0.65in}	(1+k(k+1))\hspace{0.01in}\alpha_{\scalebox{0.8}{$\displaystyle n$}}\leq \alpha_{\scalebox{0.8}{$\displaystyle n+1$}}\hspace{0.85in}\forall\hspace{0.025in} n\in I_{\scalebox{0.8}{$\displaystyle k$}},\hspace{0.05in} k\in \mathbb{N}\hspace{0.05in}
\end{equation}
\end{itemize}
The similiar observation can be given by following the same step for the case $n+1\in I_{1}$. Then, we have a regularity condition for a K\"{o}the space $\bf \mathcal{K}_{\alpha}$:
\begin{prop}\label{regularityprop} Let $\bf \mathcal{K}_{\alpha}$ be an element of the family $\bf \mathcal{K}$ parameterized by the sequence $\alpha$. Then, $\bf \mathcal{K}_{\alpha}$ is regular if and only if the inequality 
$$\hspace{0.65in}	(1+s(s+1))\hspace{0.01in}\alpha_{\scalebox{0.8}{$\displaystyle n$}}\leq \alpha_{\scalebox{0.8}{$\displaystyle n+1$}}\hspace{0.85in}$$
is satisfied for all $n\in I_{\displaystyle \scalebox{1}{$s$}}$ and $s\in \mathbb{N}$.
\end{prop}
~\\$\vspace{-.35in}$~\\
We also note that the sequence $\displaystyle \left(\alpha_{n}\right)_{n\in \mathbb{N}}=\left( \prod^{\displaystyle \scalebox{0.75}{$n-1$}}_{\displaystyle \scalebox{0.75}{$i=0$}}(1+i(i+1))\right)_{n\in \mathbb{N}}$ satisfies the condition of Proposition \ref{regularityprop} since 
$$\displaystyle {\alpha_{\displaystyle \scalebox{0.75}{$n+1$}}\over \alpha_{\displaystyle \scalebox{0.75}{$n$}}}=(1+n(n+1))\geq (1+s(s+1))$$
for all $n\in I_{\displaystyle \scalebox{0.9}{$s$}}$, $s\in \mathbb{N}$.  \par 
As a consequence of Proposition \ref{regularityprop}, we obtain the following result:
\begin{cor}\label{regularunstable} Let $\bf \mathcal{K}_{\alpha}$ be an element of the family $\bf \mathcal{K}$ parameterized by the sequence $\alpha$. If $\bf \mathcal{K}_{\alpha}$ is regular, then the sequence $\alpha$ is unstable.
\end{cor}
\begin{proof} Let $\bf \mathcal{K}_{\alpha}$ be a regular K\"{o}the space  generated by the matrix $(a_{\displaystyle \scalebox{0.75}{$k,n$}})_{k,n\in \mathbb{N}}$ given in \ref{es} and assume $\alpha$ is not unstable, that is, $\displaystyle \lim_{n\rightarrow \infty} \hspace{0.025in}{\alpha_{\displaystyle \scalebox{0.75}{$n+1$}}\over \alpha_{\displaystyle \scalebox{0.75}{$n$}}}\neq +\infty$. Then, there exist a $M>0$ and a non-decreasing sequence $\left( n_{\displaystyle \scalebox{0.75}{$k$}}\right)_{k\in \mathbb{N}}$ so that $\displaystyle  \sup_{k\in \mathbb{N}}\hspace{0.025in} {\alpha_{\displaystyle \scalebox{0.8}{$n_{\displaystyle \scalebox{0.8}{$k$}}+1$}}\over \alpha_{\displaystyle \scalebox{0.8}{$n_{\displaystyle \scalebox{0.8}{$k$}}$}}} <M.$ Since $\left( n_{\displaystyle \scalebox{0.75}{$k$}}\right)_{k\in \mathbb{N}}$ is non-decrasing and $\bf \mathcal{K}_{\alpha}$ is regular, we can write 
$$ {\alpha_{\displaystyle \scalebox{0.75}{$k+1$}} \over  \alpha_{\displaystyle \scalebox{0.75}{$k$}}} \leq {\alpha_{\displaystyle \scalebox{0.8}{$n_{\displaystyle \scalebox{0.8}{$k$}}+1$}}\over \alpha_{\displaystyle \scalebox{0.8}{$n_{\displaystyle \scalebox{0.8}{$k$}}$}}}\leq M$$
for all $k\in \mathbb{N}$ and from Proposition \ref{regularityprop}, we find that
$$\displaystyle (1+s(s+1))\leq{\alpha_{\displaystyle \scalebox{0.75}{$k+1$}} \over  \alpha_{\displaystyle \scalebox{0.75}{$k$}}} \leq M$$
for all $k\in I_{\displaystyle \scalebox{0.95}{$s$}}$, $s\in \mathbb{N}$. This is a contradiction, therefore $\alpha$ must be unstable, as desired.
\end{proof}
\begin{remark} Being unstable is not sufficient for regularity of K\"{o}the space $\bf \mathcal{K}_{\alpha}$. For instance, the sequence $\left( \alpha_{\displaystyle \scalebox{0.8}{$n$}}\right)_{n\in \mathbb{N}} =\left( n!\right)_{n\in \mathbb{N}}$ does not satisfy the condition of Proposition \ref{regularityprop}. Indeed, for every $s\in \mathbb{N}$, $\displaystyle n={s(s+1)\over 2 }\in I_{\displaystyle \scalebox{0.95}{$s$}}$ and
$$\displaystyle {\alpha_{\displaystyle \scalebox{0.8}{$n+1$}}\over \alpha_{\displaystyle \scalebox{0.8}{$n$}}}=n+1=1+{s(s+1)\over 2}< 1+s(s+1).$$
\end{remark}
\begin{remark} As a corollary of Proposition \ref{EDD}, Proposition \ref{EADD} and Corollary \ref{regularunstable}, we can obtain that $\displaystyle\Delta( \mathcal{K}_{\alpha} )=\Delta(\Lambda_{1}\left(\alpha_{\displaystyle \scalebox{0.8}{$n+1$}}\right) )$ and $\displaystyle\delta(  \mathcal{K}_{\alpha})\neq \delta(\Lambda_{1} \left(\alpha_{\displaystyle \scalebox{0.8}{$n+1$}}\right) )$ for a regular element $\bf \mathcal{K}_{\alpha}$ of the family $\bf \mathcal{K}$ which is parameterized by an exponent sequence $\alpha$.\end{remark}
\section{Some Results Obtained with the Family $\bf \mathcal{K}$}
\par In this section, we compile some additional information for the family $\bf \mathcal{K}$. We have shown that an element $\bf \mathcal{K}_{\alpha}$ of the family $\bf \mathcal{K}$ which is parametrized by an unstable sequence $\alpha$ constitutes a counterexample to Question \ref{qi2}. An element $\bf \mathcal{K}_{\alpha}$ of the family $\bf \mathcal{K}$ which is parametrized by an unstable sequence $\alpha$ is crucial for Question \ref{qi1}, as well:
\begin{thm}\label{t3} There exists a nuclear Fr\'echet space $E$ with the properties $\underline{DN}$ and $\Omega$ satisfying $\Delta(E)=\Delta(\Lambda_{1}(\varepsilon))$, for its associated exponent sequence $\varepsilon$, with the property that there is no subspace of $E$ which is isomorphic to $\Lambda_{1}(\varepsilon)$.
\end{thm}
\begin{proof} Let $\bf \mathcal{K}_{\alpha}$ be an element of the family $\bf \mathcal{K}$ which is parametrized by an unstable sequence $\alpha$. We proved that $\Delta( \mathcal{K}_{\alpha})=\Delta(\Lambda_{1}(\alpha_{\displaystyle \scalebox{0.8}{$n+1$}}))$ in Proposition \ref{EDD}. Therefore, the sequence $\left( \alpha_{\displaystyle \scalebox{0.8}{$n+1$}}\right)_{n\in \mathbb{N}}$ is the associated exponent sequence of $\bf \mathcal{K}_{\alpha}$. 
Assume that there exists a subspace of $\bf \mathcal{K}_{\alpha}$ which is isomorphic to $\Lambda_{1}(\alpha_{\displaystyle \scalebox{0.8}{$n+1$}})$. This gives us that $\delta(\Lambda_{1}(\alpha_{\displaystyle  \scalebox{0.8}{$n+1$}})) \subseteq \delta(\bf \mathcal{K}_{\alpha})$ by Proposition \ref{inv}. Since always  $\delta( \mathcal{K}_{\alpha})\subseteq \delta(\Lambda_{1}(\alpha_{\displaystyle \scalebox{0.8}{$n+1$}}))$, we conclude that $\delta( \mathcal{K}_{\alpha})= \delta(\Lambda_{1}(\alpha_{\displaystyle \scalebox{0.8}{$n+1$}}))$. But this is a contradiction since we showed that $\delta( \mathcal{K}_{\alpha})\neq \delta(\Lambda_{1}(\alpha_{\displaystyle \scalebox{0.8}{$n+1$}}))$ in Proposition \ref{EADD}. Hence, there is no subspace of $\bf \mathcal{K}_{\alpha}$ which is isomorphic to $\Lambda_{1}(\alpha_{\displaystyle \scalebox{0.8}{$n+1$}})$.
\end{proof}
\begin{remark} The above theorem indicates that Question \ref{qi1} has a negative answer. It is worth mentioning that we can find even a nuclear regular K\"{o}the space with the properties listed in Theorem \ref{t3}.
\end{remark} 
~\par In \cite{ND}, we gave conditions confirming an affirmative answer for Question \ref{qi2}. First result was related to the topology on diametral dimension of a nuclear Fr\'echet space. The diametral dimension
\begin{equation*}
\begin{split}
\Delta\left(E\right)=& \left\{\left(t_{n}\right)_{n\in \mathbb{N}}: \forall\hspace{0.05in}p\in \mathbb{N} \hspace{0.05in} \exists q>p \hspace{0.05in}\lim_{n\rightarrow\infty} t_{n}d_{n}\left(U_{q},U_{p}\right)=0\right\}\\ =& \bigcap_{p\in \mathbb{N}}\hspace{0.05in}\bigcup_{q>p}\Delta\left( U_{q}, U_{p}\right) 
\end{split}
\end{equation*}
is the projective limit of inductive limits of Banach spaces $\displaystyle \Delta\left( U_{q}, U_{p}\right)$ with the norm $\displaystyle \left\|\left(t_{n} \right)_{n}  \right\| =\sup_{n\in \mathbb{N}}\left|t_{n} \right|d_{n}(U_{q},U_{p}) $. Hence $\Delta (E)$ is a topological vector space with respect to that topology which will be called  \textit{the canonical topology}. 
\begin{thm} Let $E$ be a nuclear Fr\'echet space with properties \underline{DN} and $\Omega$ and $\varepsilon=\left(\varepsilon_{n}\right)_{n\in \mathbb{N}}$ be the associated exponent sequence of $E$. If $\Delta\left(E\right)$, with the canonical topology, is barrelled, then $\Delta\left(E\right)=\Delta\left(\Lambda_{1}\left(\varepsilon\right)\right)$ if and only if $\delta\left(E\right)=\delta\left(\Lambda_{1}\left(\varepsilon\right)\right)$.
\end{thm} 
\begin{proof} \cite[Theorem 4.2]{ND}
\end{proof}
Therefore, we obtain the following:
\begin{prop} Let $\bf \mathcal{K}_{\alpha}$ be an element of the family $\bf \mathcal{K}$ parameterized by an unstable sequence $\alpha$. Then $\Delta( \bf \mathcal{K}_{\alpha})$, with the canonical topology, is neither barrelled nor ultrabornological.
\end{prop}
We actually wanted the barrelledness in \cite[Theorem 4.2]{ND} to be able to use a closed graph type theorem, \cite[Theorem 5, Pg. 40]{T} which says that a linear map $f$ from a barrelled space $X$ into a Fr\'echet space $Y$ is continuous provided that the graph of f is closed in $X\times Y$. Since $\delta( \mathcal{K}_{\alpha})\neq \delta(\Lambda_{1}(\alpha_{\displaystyle \scalebox{0.8}{$n+1$}}))$ and $\Delta( \mathcal{K}_{\alpha})=\Delta(\Lambda_{1}(\alpha_{\displaystyle \scalebox{0.8}{$n+1$}}))$, the technique used in the proof of \cite[Theorem 4.2]{ND} is not valid for an element $\bf \mathcal{K}_{\alpha}$ of the family $\bf \mathcal{K}$ parameterized by an unstable sequence $\alpha$. Hence, this gives us that the identity mapping from $\Delta(\bf \mathcal{K}_{\alpha})$ into $\Lambda_{1}(\alpha_{\displaystyle \scalebox{0.8}{$n+1$}})$ is not continuous although it has a closed graph: 
\begin{thm} Let $\bf \mathcal{K}_{\alpha}$ be an element of the family $\bf \mathcal{K}$  parameterized by an unstable sequence $\alpha$. Then  $\Delta( \mathcal{K}_{\alpha})=\Delta(\Lambda_{1}(\alpha_{\displaystyle \scalebox{0.8}{$n+1$}}))$ and the identity map from $\Delta(\bf \mathcal{K}_{\alpha})$ into $\Lambda_{1}(\alpha_{\displaystyle \scalebox{0.8}{$n+1$}})$ is not continuous although it has a closed graph. 
\end{thm}
In \cite{T2}, T. Terzio\u{g}lu defined the notion \textit{prominent bounded subset} in order to show that the diametral dimension of some Fr\'echet spaces is determined by a single bounded set:
\begin{defnt} Let $E$ be a Fr\'echet space. A bounded set $B$ is said to \textbf{prominent} if 
$$\displaystyle \Delta(E)=\left\lbrace\hspace{0.025in} \left(x_{n}\right)_{n\in \mathbb{N}}\hspace{0.05in}:\hspace{0.05in} \lim_{n\rightarrow +\infty} x_{n}\hspace{0.025in}d_{n}\left(B, U_{p}\right)=0\hspace{0.1in} \forall \hspace{0.025in}p \right\rbrace. $$
\end{defnt}
The existence of a prominent
bounded subset in the nuclear Fr\'echet space E plays a decisive role for the affirmative answer of Question \ref{qi2}.
\begin{thm}\label{t5} Let $E$ be a nuclear Fr\'echet space with the properties \underline{DN} and $\Omega$ and $\varepsilon$ the associated exponent sequence.  $\delta\left(E\right)=\delta\left(\Lambda_{1}\left(\varepsilon\right)\right)$ if and only if $E$ has a prominent bounded set  and $\Delta\left(E\right)=\Lambda_{1}\left(\varepsilon\right)$.
\end{thm}
\begin{proof} \cite[Theorem 4.8]{ND}
\end{proof}
Obviously, this condition is not valid for an element $\bf \mathcal{K}_{\alpha}$ of the family $\bf \mathcal{K}$ which is parameterized by an unstable sequence $\alpha$ since $\Delta( \mathcal{K}_{\alpha})=\Delta(\Lambda_{1}(\alpha_{\displaystyle \scalebox{0.8}{$n+1$}})$ and  $\delta( \mathcal{K}_{\alpha})\neq \delta(\Lambda_{1}(\alpha_{\displaystyle \scalebox{0.8}{$n+1$}}))$.
\begin{thm} There exists a nuclear Fr\'echet space $E$ with the properties $\underline{DN}$ and $\Omega$ satisfying $\Delta(E)=\Delta(\Lambda_{1}(\varepsilon))$ for its associated exponent sequence $\varepsilon$ such that there is no prominent bounded set of $E$.
\end{thm}
\begin{remark}\label{stabled2} It is worth to note that as a consequence of Theorem \ref{t5} and Corollary \ref{sc}, an element $\bf \mathcal{K}_{\alpha}$ of the family $\bf \mathcal{K}$ parameterized by a \underline{stable} sequence $\alpha$ has a prominent bounded subset.
\end{remark}
\par A nuclear Fr\'echet space $E$ with an increasing sequence  of seminorms $\left(  \left\| .\right\|_{k} \right)_{k\in \mathbb{N}}$ is called \textit{tame} if there exists an increasing function $\sigma:\mathbb{N}\rightarrow \mathbb{N}$, such that for every continuous linear operator $T:E\rightarrow E$ there exists a $n_{0}\in \mathbb{N}$ and $C>0$ so that
$$\hspace{1.65in}\left\|T(x) \right\|_{k}\leq C\left\| x\right\|_{\sigma \left( k\right) } \hspace{1.5in} \forall \hspace{0.045in}x\in E.$$
\par In \cite[Theorem 2.3]{A}, A. Aytuna proved that a nuclear Fr\'echet space $E$ with the properties $\underline{DN}$ and $\Omega$ and stable associated exponent sequence $\varepsilon$ is isomorphic to a power series space of finite type if and only if $E$ is tame and $\delta(E)=\delta(\Lambda_{1}(\varepsilon))$. As a consequence of this result and Remark \ref{stabled2}, we have the following: 
\begin{prop}  Let $\bf \mathcal{K}_{\alpha}$ be an element of the family $\bf \mathcal{K}$ parameterized by a stable sequence $\alpha$. Then, $\bf \mathcal{K}_{\alpha}$ is not tame.
\end{prop}


\end{document}